\documentclass[12pt,reqno]{amsart}
\usepackage{VMpreamble}

\usepackage{amsmath,amssymb,amsfonts,amscd,mathtools}
\usepackage{fullpage}
\usepackage[indent]{parskip}
\usepackage{comment}
\usepackage{ytableau}
\usepackage{array}
\definecolor{almond}{rgb}{0.94, 0.87, 0.8}
\definecolor{crimson}{rgb}{0.86, 0.08, 0.24}
\definecolor{mintcream}{rgb}{0.96, 1.0, 0.98}
\definecolor{champagne}{rgb}{0.97, 0.91, 0.81}
\definecolor{bubblegum}{rgb}{0.99, 0.76, 0.8}
\definecolor{airforceblue}{rgb}{0.36, 0.54, 0.66}
    	\definecolor{almond}{rgb}{0.94, 0.87, 0.8}
    	\definecolor{antiquefuchsia}{rgb}{0.57, 0.36, 0.51}
    	\definecolor{antiquewhite}{rgb}{0.98, 0.92, 0.84}

\newcommand{\ev}[1]{\Big\langle #1 \Big\rangle}

\newcommand\lam{\lambda}

\newcommand\bx{\mathbf{x}}
\newcommand\by{\mathbf{y}}

\newcommand\sfZ{\mathsf{Z}}

\newcommand\sfC{\mathsf{C}}

\newcommand{\bfx}{\mathbf{x}}

\newcommand{\bfy}{\mathbf{y}}

 \newcommand{\Rem}[1]{\operatorname{Rem}(#1)}
 \newcommand{\Add}[1]{\operatorname{Add}(#1)}

\usetikzlibrary{arrows.meta}

\definecolor{FusionInk}{HTML}{17202A}
\definecolor{FusionBlue}{HTML}{24445C}
\definecolor{FusionArrow}{HTML}{59636B}

\newtheorem{definition}{Definition}[section]
\newtheorem*{definition*}{Definition}
\newtheorem{proposition}[definition]{Proposition}

\newtheorem{theorem}[definition]{Theorem}

\newtheorem*{theorem*}{Theorem}

\newtheorem{corollary}[definition]{Corollary}
\newtheorem{lemma}[definition]{Lemma}

\newcommand{\Wd}{W^{\rm d}}
\newcommand{\Wf}{W^{\rm f}}
\newcommand{\Wp}{W^{\rm p}}

\title[]{Bases for shift invariant symmetric polynomials}

\author{Victor Mishnyakov}
\address[V.M.]{Moscow Institute of Physics and Technology,
Dolgoprudny 141701, Russia \& NRC ``Kurchatov Institute'',
Moscow 123182, Russia}
\email{mishnyakovvv@gmail.com}

\author{Nikita Tselousov}
\address[N.T.]{Moscow Institute of Physics and Technology,
Dolgoprudny 141701, Russia;  Institute for Information Transmission Problems
of the Russian Academy of Sciences, Moscow 127051, Russia \& NRC ``Kurchatov Institute'',
Moscow 123182, Russia}
\email{tselousov.ns@phystech.edu}

\newcommand{\preprintnumber}{%
  \shortstack[r]{MIPT/TH-31/26 \\ IITP/TH-34/26\\ ITEP/TH-40/26 }}
  \makeatletter
\let\original@settitle\@settitle
\renewcommand{\@settitle}{%
  \noindent\hfill{\normalfont\small\preprintnumber}\par
  \vspace{1em}%
  \original@settitle
}
\makeatother

\makeatletter

\renewcommand{\l@section}{%
  \@tocline{1}{2pt}{0pt}{}{\small}}

\renewcommand{\l@subsection}{%
  \@tocline{2}{0pt}{1.5em}{}{\footnotesize}}

\makeatother

\begin{document}
\raggedbottom
\begin{abstract}
We revisit a special class of symmetric polynomials of $N$ variables that are invariant under simultaneous shifts $x_i \to x_i + c$ of all variables. Our main goal is to look at the Schur expansion of such symmetric functions. We construct new bases in the space of shift invariant symmetric polynomials and revisit some previously known ones. We call them the fusion, fission, decay, and character bases. The basis elements have explicit and computable expansion coefficients in the Schur polynomial basis. The character and decay bases are only defined for polynomials of degree $d \leq N$. On the other hand, the fusion and fission bases exist in every degree. Moreover, we find that after a proper rescaling, they define a family of symmetric functions, depending on an extra parameter $u$ that interpolates between all $N$. Shift invariant symmetric polynomials have appeared in various contexts, in particular, in problems motivated by physics. As such, we briefly discuss several properties and applications of our bases,  relevant in different scenarios.

\end{abstract}

\maketitle
\begingroup
  \setlength{\parskip}{0pt}
  \linespread{1}\selectfont
  \tableofcontents
\endgroup

\section{Introduction}
In this paper we study the space of symmetric polynomials that are invariant under a simultaneous constant shift of all their variables. Such symmetric polynomials appeared in various contexts in mathematics and physics. In particular, they appear in the study of systems of interacting particles in the center-of-mass frame, in systems like the Calogero model and in certain formulations of the fractional quantum Hall effect (QHE). Within the $M$-body generalization of Haldane's pseudopotential approach one is required to identify a basis of states for the relative motion part of the Hilbert space, which for bosons amounts exactly to finding a basis in the space of shift invariant symmetric polynomials \cite{simon2007pseudopotentials}. Another natural occurrence of such polynomials, related to the fractional QHE, is in wavefunctions given by Jack polynomials at special negative rational values of $\beta$ and indexed by partitions of special shape \cite{bernevig2008model}. Within these contexts other interesting properties such as clustering and squeezing are discussed \cite{liptrap2008translation,baratta2011jack}. Clustering properties of Jack polynomials at these special values were considered in \cite{feigin2002differential}. Shift invariant symmetric polynomials appear in the description of the center of universal enveloping algebras of special linear groups where tug-the-hook symmetry effectively acts as translations in terms of Frobenius hook-variables \cite{Mishnyakov:2020khb}, a construction which plays a role in the theory of quantum knot invariants. The same question appears when one expands BPS operators in an $SU(N)$ supersymmetric Yang-Mills theory in terms of their $U(N)$ counterparts \cite{deMelloKoch:2004crq}. Special bases in the space of shift invariant symmetric polynomials also originate from symmetrized graph monomials \cite{alexandersson2012discriminants}.

Our initial motivation to study shift invariant symmetric polynomials comes from the construction utilized in \cite{Dixon:2026fuz}. In that paper certain combinations of classical polylogs, called ladder integrals, were studied. It was observed that the antipode map of the multiple polylog Hopf algebra, restricted to single polylogs, can be mapped to the shift operator on symmetric polynomials. Hence every shift invariant symmetric polynomial corresponds to an antipodally self dual polynomial in the ladder integrals. However, a general statement was found  only for those symmetric polynomials where the number of variables $N$ is greater or equal than the total degree of the polynomials. The results of the present paper complete the construction. First, we construct another basis that is connected to forgotten symmetric functions. Second, we use the involution in the space of symmetric functions to construct a dual pair for a basis considered in \cite{Dixon:2026fuz}. Both these constructions are applicable for arbitrary $N$.

Let $\Lambda_N$ be the ring of symmetric polynomials in $N$ variables $\bx=(x_1,\ldots, x_N)$ and $ {\bf 1} = (1,\ldots, 1)$. Then denote the ring of shift invariant symmetric polynomials in $N$ variables as $\bar{\Lambda}_{N}$:
\begin{equation}
    \bar{\Lambda}_{N} = \left\{ f(\bx) \in \Lambda_N :  f(\bfx+c {\bf 1})=f(\bfx)  \,,\forall c \in \mathbb{C} \right\}
\end{equation}
Just as with usual symmetric polynomials, this space is naturally graded by the degree of the polynomial:
\begin{equation}
    \bar{\Lambda}_{N} = \bigoplus_{d=0}^{\infty}\bar{\Lambda}_{N,d} 
\end{equation}

Linearly independent shift invariant symmetric polynomials are labeled by \emph{nuclear} partitions \cite{schneider2019nuclear}, which are partitions with components greater than one. For $N$ variables the components are also at most $N$ \footnote{This statement is similar to the well known fact in representation theory of $GL(N)$ where characters -- Schur polynomials -- are enumerated by partitions of restricted length $l(\lambda) \leq N$.}, but without restrictions on the length of the partition. 
Denote the set of all nuclear partitions of $d$ as:
\begin{equation}
    \mathcal{N}(d):= \left\{ \lambda \vdash d : \lambda_i \geq 2 \right\}
\end{equation}
and those with all row length greater than one and less than or equal to $N$ as:
\begin{equation}
    \mathcal{N}_N(d):= \left\{ \lambda \vdash d : 2 \leq\lambda_i \leq N \right\}
\end{equation}
Set $\mathcal{N}(0)=\{\varnothing\}$.
There are several natural ways to construct a basis in $\bar{\Lambda}_N$ using center-of-mass coordinates \cite{liptrap2008translation}, one of which we will review below. The main goal of this paper is to construct bases that are explicit as functions of the initial variables $\bx$. That is, we would like to find such bases whose expansion coefficients into some well-understood natural basis in the space of symmetric polynomials would be known explicitly. In particular, we choose the basis of Schur polynomials $s_{\lambda}(\bfx)$. Therefore we want to find, for every $\lambda$ in $\mathcal{N}_N(d)$, a shift invariant symmetric polynomial
\begin{equation}
    W_{\lambda}(\bx)=\sum_{|\mu|=|\lambda|} \alpha_{\lambda\mu}s_{\mu}(\bfx)
\end{equation}
and to find explicit expressions for coefficients $\alpha_{\lambda\mu}$. 

In this paper we present four such possibilities each distinguished by certain special properties. To define them let us introduce a notation:
\begin{equation}\label{eq:C}
    \sfC_{\lambda}(N)=\prod_{\Box \in \lambda}(N+j_{\Box}-i_{\Box})
\end{equation}
Here $i_{\Box}$, $j_{\Box}$ are the row and column numbers of the box, respectively, in the Young diagram. The coordinates of the first corner box in any Young diagram are $(i_{\Box}, j_{\Box}) = (1,1)$.
 
Our main result concerns the construction of the  \emph{fusion basis}, whose elements are the symmetric polynomials defined by:
\begin{equation}
        \widetilde{W}_{\lambda}^{\rm f}= D_{\lambda}(N) \sum_{\substack{\mu  \vdash|\lambda|\\
        \ell(\mu)\leq N}} K^{-1}_{\lambda\mu^T} \, \dfrac{ s_\mu(\bx)} {\sfC_{\mu}(N)}
    \end{equation}
where $\mu^T$ is the transpose partition, $K^{-1}_{\lambda\mu}$ are elements of the inverse Kostka matrix and $D_{\lambda}(N)$ is a common factor defined in sec.\ref{sec::uniform N}.

Another basis in the space of symmetric functions comes from the decay and fission constructions. An alternative version of it was first proposed in \cite{deMelloKoch:2004crq} without a general proof of shift invariance. Here we give an explicit combinatorial formula and prove its shifts invariants. We will refer to this family of polynomials as the \emph{fission basis}. Denote by $\eta^{(k)}$ the partition obtained from $\eta$ as $\eta^{(k)}=(\eta_1+k-|\eta|,\eta_2,\ldots )$, $|\eta|\leq k$. Then, define:
\begin{equation}
    \label{fission def}
\widetilde{W}^{\rm F}_{\lambda}(\mathbf x)
= D^{\prime}_{\lam} (N) \sum_{\substack{\eta\subseteq\lambda^T\\
                  \eta_1=(\lambda^T)_1}}
(-1)^{d-|\eta|}f^{\lambda^T/\eta} \,
\frac{s_{\eta^{(d)}}(\mathbf x)}
     {\sfC_{\eta^{(d)}}(N)},
\qquad d=|\lambda|
\,,
\end{equation} 
where $f^{\lambda/\nu}$ is the number of standard Young tableaux \footnote{This notation should not be confused with forgotten symmetric functions that we introduce later in sec.\ref{sec::forgotten}.} of skew shape $\lambda/\nu$.

\begin{theorem}\label{thrm:fission}
    Each of the two families of symmetric polynomials $\widetilde{W}_{\lambda}^{\rm f}(\bx)$ and $\widetilde{W}_{\lambda}^{\rm F}(\bx)$ for $\lambda \in \mathcal{N}_N$ is a basis of the space of shift invariant symmetric polynomials $\bar{\Lambda}_N$. 
\end{theorem}
Denote by $\Lambda_{N,\leq d}$ the subspace:
\begin{equation}
    \Lambda_{N,\leq d}=\bigoplus\limits_{j=0}^{d} \Lambda_{N,j}
\end{equation}
and $\bar{\Lambda}_{N,\leq d}$ its shift invariant subspace. For symmetric polynomials with degree less than or equal to $N$ there exist two more distinguished choices of basis.
\begin{itemize}
    \item The \emph{decay basis}. Denote by $\nu^{[k]}$ the partition obtained from $\nu$ by appending  $(k-|\nu|)$ ones (i.e. rows of length one) if $|\nu| \leq k$. 
    \begin{equation}
        W^{\rm d}_{\lambda}(\bx)= \sum_{\substack{\nu\subseteq \lambda\\
        \ell(\nu)=\ell(\lambda)}} (-1)^{|\lambda|-|\nu|} \, f^{\lambda/\nu} \, \dfrac{s_{\nu^{[|\lambda|]}}(\bx)}{\sfC_{\nu^{[|\lambda|]} }(N)} 
    \end{equation}
    \item The \emph{character basis}:
    \begin{equation}
        W_{\lambda}^{\rm p}(\bx) = \sum_{\mu  \vdash|\lambda|}\chi_{\lambda}^{\mu} \, \dfrac{s_{\mu}(\bx)}{\sfC_\mu(N)}
    \end{equation}
    where $\chi_{\lambda}^{\mu}$ are symmetric group characters in representation $\mu$, corresponding to cycle type $\lambda$.
\end{itemize}
For these we will prove the following.
\begin{theorem}\label{thrm:decay}
 For $0\leq d\leq N$ the families $ \Wd_{\lambda}(\bx)$ and $ \Wp_{\lambda}(\bx)$  for $\lambda \in  \bigcup_{s=0}^{d}\mathcal{N}(s)$ are each bases in the space $\bar{\Lambda}_{N,\leq d}$.
\end{theorem}
For both bases the result has been obtained earlier. The decay basis part of the theorem was already proved in \cite{Dixon:2026fuz}, where the decay basis was first discovered, without its relation to the fission basis. For the character basis it was proved in \cite{Brown:2007bb}. Here we provide a new proof for the character basis that is uniform with all the other families. We will see, that in a certain precise sense the \textit{fission basis is dual to the decay basis} under the involution on symmetric functions. Different bases that we constructed here are distinguished by different features. In particular the coefficients of the decay and fission bases $f^{\lambda/\nu}$ have a closed form determinant formula and therefore can be easily computed. The Kostka numbers and characters of the symmetric group that appear in the fusion and character bases are computationally more involved. On the other hand as we will see  the fusion and character bases are very natural from the point of view of our construction.

\section{Symmetric functions and shift invariance}
In this section we introduce standard notations and introduce shift invariant symmetric polynomials and functions. Let $\Lambda_N=\mathbb{C}[x_1,\ldots x_N]^{S_N}$ denote the ring of symmetric polynomials in $N$ variables and let $\Lambda$ be the stable ring of symmetric functions. Denote the projection as $\pi_N:\Lambda \rightarrow \Lambda_N$.
\\
The ring $\Lambda_N$ is graded by degree of the polynomial. Denote the respective graded component by $\Lambda_{N,d}$.

Let $\mathcal{P}(d)$ denote the set of integer partitions of $d$ and $\mathcal{P}=\bigcup_{d=0}^{\infty} \mathcal{P}(d)$. Similarly, let $\mathcal{P}_N(d)$ denote the set of partitions with components less than or equal to $N$:
\begin{equation}
    \mathcal{P}_N(d) = \left\{ \lambda \vdash d: \lambda_i \leq N\right\}
\end{equation}
Denote the Hall inner product \cite{Macdonald:book} on the space of symmetric functions as $(f,g)$. It is defined by requiring orthogonality of the power sums:
\begin{equation}
    (p_{\lambda},p_\mu)=z_{\lambda} \delta_{\lambda\mu}
\end{equation}
where $z_{\lambda} = \prod_i m_{i}(\lambda)! i^{m_{i}(\lambda)}$ is a standard factor and $m_i(\lambda)$ is the multiplicity of part $i$ in the partition $\lambda$.   
Among symmetric functions (and polynomials) we are looking for those that obey an extra condition:
\begin{equation}
    f(\bx+c{\bf 1})=f(\bx)
\end{equation}
where by $\bx+c$ we mean a simultaneous shift of all variables, i.e. $x_i+c$. Before going to more general aspects, let us look at some examples. Obvious examples of such polynomials are even powers of the Vandermonde determinant:
\begin{equation}
    \Delta^{2n}(\bx)=\prod_{1\leq i<j\leq N} (x_i-x_j)^{2n}
\end{equation}
For the case of two variable polynomials the whole space $\bar{\Lambda}_2$ is spanned by even powers of the Vandermonde determinant. Therefore:
\begin{equation}
    \bar{\Lambda}_{2,d}=\left\{\begin{array}{cc}
\{0\} , \quad &d \text{ -- odd}
    \\    
    \operatorname{span}\{\Delta^d(\bx)\} , \quad & d \text{ -- even}
    \end{array} \right.
\end{equation} 
Going in an orthogonal direction, staying in degree two but for generic $N$, i.e. in $\bar{\Lambda}_{N,2}$, we have simple example:
\begin{equation}
   \sum_{1\leq i<j\leq N} (x_i-x_j)^2
\end{equation}
which is shift invariant by construction. One could also approach this problem in a more systematic way. Namely, since the Schur polynomials form a basis of symmetric polynomials, a generic polynomial at degree two would be a combination of $s_{[2]}(\bx)$ and $s_{[1,1]}(\bx)$. The condition for shift invariance then reads
\begin{equation}
    \alpha_{[2]} s_{[2]}(\bx)+\alpha_{[1,1]} s_{[1,1]}(\bx)=\alpha_{[2]} s_{[2]}(\bx+c{\bf 1})+\alpha_{[1,1]} s_{[1,1]}(\bx+c{\bf 1})\,, \quad \forall c
\end{equation}
One can check that even though the equation imposes an equality in both degrees zero and one, there is only one linearly independent constraint. This condition can be solved and the coefficients are found up to an overall normalization. For $N \geq 2$ it can be chosen such that:
\begin{equation}
    W_{[2]}(\bx)=\dfrac{ s_{[2]}(\bx)}{N(N+1)}-\dfrac{ s_{[1,1]}(\bx)}{N(N-1)}
\end{equation}
In terms of squares of differences this gives:
\begin{equation}
      W_{[2]}(\bx)=\dfrac{ \sum\limits_{1\leq i<j\leq N} (x_i-x_j)^2}{N(N^2-1)}
\end{equation}
To understand the dimension of the homogeneous components $\bar{\Lambda}_{N,d}$ one can use center-of-mass coordinates (COM). This idea was developed in \cite{liptrap2008translation}. Define:  
\begin{equation}\label{eq:y}
    y_i=x_i- \dfrac{1}{N}\sum_{j=1}^N x_j=x_i- \bar{x}\,.
\end{equation}
Let $ \pi^{\rm COM}$ be the corresponding ring homomorphism:
\begin{equation}
\begin{gathered}
    \pi^{\rm COM}: \     \Lambda_N \rightarrow \Lambda_N
    \\
    f(\bx) \mapsto f(\bx- \bar{x} {\bf 1})=f(\by)
\end{gathered}
\end{equation}

\begin{proposition}
The image of $\pi^{\rm COM}$ is the space of shift invariant symmetric polynomials:
\begin{equation}
    \operatorname{im} \pi^{\rm COM}= \bar{\Lambda}_N
\end{equation}
In other words, any symmetric polynomial in the variables $\by=\{y_1,\ldots y_N\}$ \eqref{eq:y} is a shift invariant polynomial in the $\bx$ variables. Conversely, any shift invariant polynomial can be written as a symmetric polynomial of $\by$.
\end{proposition}
\begin{proof}
    Symmetric polynomials of the $\by$ variables are shift invariant in $\bx$ by construction. Each variable $y_i$ remains unchanged under simultaneous shifts of $x_i$'s. Conversely, let $f(\bx)$ be a shift invariant polynomial, i.e., $f(\bx+c{\bf 1})=f(\bx)$. Choosing $c=-\bar{x}$, we get $f(\bx)=f(\bx-\bar{x}{\bf 1})=f(\bfy)$.
\end{proof}
Therefore center-of-mass coordinates naturally parametrize the space of shift invariant symmetric polynomials. In order to construct a basis we need to find the kernel of this map. Clearly $p_1(\bfy)=0$ by definition. The following statement was proven in 
 \cite{liptrap2008translation}.
\begin{lemma}\label{th:COMker}
    The kernel of the $\rm COM$ map is the ideal generated by  $p_1(\bx)$.
    \begin{equation}
        \ker \pi^{\rm COM} = (p_1(\bx) )= p_1(\bx)\Lambda_N
    \end{equation}
\end{lemma}
\begin{proof}
    The proof in \cite{liptrap2008translation} adapted to our notation is the following. Let $q_k=p_k(\bx-\bar{x}{\bf1})$. This change of variables is triangular with diagonal 1. Hence:
    \begin{equation}
        \Lambda_N=\mathbb{C}[p_1,q_2,\ldots, q_N]
    \end{equation}
    Under the COM map we have:
    \begin{equation}
        \pi^{\rm COM}(p_1)=0 \,, \quad \pi^{\rm  COM}(q_k)=q_k
    \end{equation}
\end{proof}
Therefore we have the following corollary
\begin{corollary}
  A basis in the ring of shift invariant symmetric polynomials is given by power sums of the center-of-mass coordinates:
   \begin{equation}\label{eq:powersumCOM}
        \bar{\Lambda}_{N} = \langle \, p_{\lambda}\left(\bx-\bar{x}{\bf 1}\right)  | \ \lambda \in \mathcal{N}_N \rangle
    \end{equation}
     which are labeled by nuclear partitions. We have the corresponding ring isomorphism:
    \begin{equation}
        \bar{\Lambda}_N \simeq \Lambda_N/(p_1)
    \end{equation}
\end{corollary}
The dimension of the homogeneous components  is given by the following Hilbert series:
\begin{equation}
    \sum_{d \geq 0} \dim \bar{\Lambda}_{N,d} t^d = \prod_{k=2}^{N} \dfrac{1}{(1-t^k)}
\end{equation}
For $d \leq N$ the dimension equals the number of nuclear partitions of $d$ without a restriction on row lengths:
    \begin{equation}
        \dim \bar{\Lambda}_{N,d}= \big|\mathcal{N}(d)\big|= p(d)-p(d-1)
    \end{equation}
where $p(n)$ is the number of partitions of $n$. Therefore in order to label basis elements for shift invariant polynomials it is enough to construct a basis for polynomials of $\by$. In particular, one can choose the power sum basis:
\begin{equation}
    p_{\lambda }(\by)=\prod_{i=1}^{\ell(\lambda)} p_{\lambda_i}(\by)\,, \quad \lambda\in \mathcal{P}
\end{equation}
Notice, however, that $p_1(\by)=0$, which means that to obtain non-zero elements the partitions must have $\lambda_i > 1$ for all $i$. In addition, the standard restriction is that $\lambda_i \leq N$ for all $i$, for the corresponding power sums to be algebraically independent. However, one must also make sure that no extra relations between $p_{\lambda }(\by)$ appear, which is guaranteed by Lemma \ref{th:COMker}. 

In principle, this center-of-mass construction gives a perfectly fine basis. Suppose, however, that we would like to understand its structure as a function of the original $x_i$ variables. For example, let us ask how it expands into the basis of Schur polynomials. Then we will find that the coefficients of this expansion do not have a nice structure. For example:
\begin{equation}
\begin{split}
        &p_4(\bx-\bar{x})=p_4-\dfrac{4p_1 p_3}{N}+\dfrac{6p_1^2 p_2}{N^2}-\dfrac{3p_1^4}{N^3}=\\&
        =\frac{\left(N^3-4 N^2+6 N-3\right) s_{[4]}(\bx)}{N^3}-\frac{\left(N^3-6 N+9\right) s_{[3,1]}(\bx)}{N^3}+\frac{2 \left(2 N^2-3\right) s_{[2,2]}(\bx)}{N^3}\\&+\frac{\left(N^3-6 N-9\right)
   s_{[2,1,1]}(\bx)}{N^3}-\frac{\left(N^3+4 N^2+6
   N+3\right) s_{[1,1,1,1]}(\bx)}{N^3}
\end{split}
\end{equation}
Although in this lower degree example the coefficients are not too complicated, we can already observe that they are expressed in terms of nonfactorizable polynomials in $N$. In principle, one could brute force such an expansion in general, by iteratively expanding the binomials $x_i-\bar{x}$. However, the resulting coefficients are extremely complicated.

This discussion motivated the question that we raise in this paper. Can we find a basis in the space of shift invariant symmetric polynomials that would have an expansion in terms of Schur polynomials, with coefficients that have a simple form? Clearly, there are arbitrary choices that were made here. The first one is the choice of Schur polynomials as the basis for the expansion. Perhaps simpler expansions exist in some of the other well-known bases. Second, we must specify what we mean by simplicity of the coefficients. The construction below will make it clear that we mean that the dependence on $N$ is fully factorized.
\\
\subsection{Construction of shift invariant symmetric polynomials}

Now let us explain the general procedure to find shift invariant symmetric polynomials with a factorized Schur expansion. We first explain the construction schematically. Let 
\begin{equation}
    \nabla_N= \sum_{i=1}^N \dfrac{\partial}{\partial x_i}
\end{equation}
be the infinitesimal shift operator. We will show that, up to conjugation by an operator $\hat{\sfC}(N)$, the shift operator is equivalent to the derivative with respect to the first power sum $\dfrac{\partial}{\partial p_1}$ in the stable ring of symmetric functions. This means that symmetric functions independent of the first power sum will produce shift invariant symmetric polynomials. At the same time $\hat{\sfC}(N)$ is diagonal on Schur functions with a factorized eigenvalue $\sfC_\lambda(N)$ from \eqref{eq:C}. This endows the constructed shift invariant functions with a simple Schur expansion, where all the dependence on $N$ enters only via the eigenvalues of $\hat{\sfC}(N)$.

Now, we proceed step by step. For $u \in \mathbb{C}$  define an operator:
\begin{equation}
    \hat{\sfC}(u): \Lambda \rightarrow \Lambda
\end{equation}
such that it is diagonal on Schur functions:
\begin{equation}
    \hat{\sfC}(u) s_{\lambda} = \sfC_{\lambda}(u) s_{\lambda}
\end{equation}
where
\begin{equation}\label{eq:Cev}
    \sfC_{\lambda}(u)=\prod_{\Box \in \lambda}(u+j_{\Box}-i_{\Box})
\end{equation}
For generic $u$ the operator has a trivial kernel. However, we will need to evaluate it at $u=N$ - an integer value. At positive integer values of $u$ the eigenvalues \eqref{eq:Cev} can vanish. This happens because when $\ell(\lambda)>N$, there exists a box of the Young diagram with coordinates $(i=N+1,j=1)$. 
\begin{proposition}
    The operator $\hat{\sfC}(N)$ has a kernel that is spanned by Schur functions corresponding to long enough partitions:
\begin{equation}
   \ker \hat{\sfC}(N) = \operatorname{span} \left\{ s_{\lambda}:\ell(\lambda) > N  \right\} = I_N
\end{equation}
\end{proposition}
The kernel $I_N$ is an ideal in $\Lambda$ due to the Littlewood-Richardson rule \cite{Macdonald:book}. The finite variable projection $\pi_N$ acts on Schur functions as \cite{Macdonald:book}:
\begin{equation}
         \pi_N (s_{\lambda}) = 
         \left\{ \begin{split}
            & s_{\lambda}(\bx)  \,,\ &\ell(\lam) \leq N
             \\
             &0 \,,\  &\ell(\lam) > N
         \end{split} \right.
    \end{equation}
This means that we also have:
\begin{equation}
    I_N = \ker \pi_N
\end{equation}
Denote by $V_N$ the subspace of $\Lambda$ that is spanned by Schur functions with partitions of length at most $N$:
\begin{equation}
    V_N=\operatorname{span} \{ s_\mu : \ell(\mu) \leq N \} \subset \Lambda
\end{equation}
\begin{proposition}
With respect to the Hall inner product \cite{Macdonald:book} we have a vector space isomorphism
\begin{equation}
    V_N=I_N^{\perp}
\end{equation}
and 
\begin{equation}
    \Lambda=V_N\oplus I_N
\end{equation}
\end{proposition}
Finally, restricting the projection to $V_N$ we get the graded vector space isomorphism \footnote{$\Lambda_N$ is also isomorphic to $\Lambda/I_N$ as  rings.}:
\begin{equation}
    \pi_N\big|_{V_N}:V_N \simeq\Lambda_N
\end{equation}
Now we come to the key statement of this section. To motivate and set it up, we first establish and compare the action of the shift operator $\nabla_N$ and the operator $ \dfrac{\partial}{\partial p_1}$ in the Schur basis. This action is convenient to formulate via the set $\Rem\lambda$ that is a subset of boxes of the diagram $\lambda$ that can be removed so that the resulting diagram remains a Young diagram. In terms of their coordinates:
\begin{equation}
    \Rem\lambda = \left\{ (i,j) \in \lambda \, : \,  j=\lambda_{i} \leq \lambda_{i-1},\lambda_i >\lambda_{i+1}  \right\}
\end{equation}
These are known classical results that can be found in \cite{Macdonald:book}.
\begin{lemma}
    The shift operator has the following expansion:
    \begin{equation}
        \nabla_N s_{\lambda}(\bx)=\sum_{\Box \in \Rem\lambda}(N+j_\Box-i_\Box) s_{\lambda-\Box}(\bx)
    \end{equation}
    and in the stable ring we have the dual Pieri rule:
\begin{equation}\label{eq:DualPieri}
    \dfrac{\partial}{\partial p_1} s_{\lambda}=\sum_{\Box \in \Rem\lambda}s_{\lambda-\Box}
\end{equation}
\end{lemma}
\begin{proof}
    It is a classical consequence of the  finite shift binomial formula for Schur polynomials (Section I.3, Example 10 in \cite{Macdonald:book}) :
 \begin{equation}
        s_\lambda(\bfx+a{\bf1}) = \sum_{\mu\,:\, \mu\subseteq \lambda} a^{|\lambda|-|\mu|}b_{\lambda\mu} s_\mu(\bfx) \,,
    \end{equation}
   where $b_{\lambda\mu}$ is given by:
    \begin{equation} 
        b_{\lambda\mu} = \det\left[ \binom{\lambda_i-i+N}{\mu_j-j+N} \right] \,.
    \end{equation}
    When expanded to first order in $a$ it leads to the expression for the shift operator. 
    \\
    
    The dual Pieri rule follows from the fact that the operator $\dfrac{\partial}{\partial p_1}$ is adjoint to multiplication by $p_1$ w.r.t. to the Hall inner product. The Schur expansion of multiplication by $p_1$ is controlled by the Pieri rule.
\end{proof}

Comparing the two expressions we see that the extra factor $(N+j_\Box-i_\Box)$ can be accounted for by the action of the $\hat{\sfC}(N)$ operator. The following lemma establishes the precise relation.
\begin{lemma}\label{lemma:transfer}
    On $V_N$ we have:
    \begin{equation}
        \nabla_N  \pi_N \hat{\sfC}(N)^{-1} \Big|_{V_N}=\pi_N  \hat{\sfC}(N)^{-1} \dfrac{\partial}{\partial p_1} \Bigg|_{V_N}
    \end{equation}
And therefore the operator $\pi_N \hat{\sfC}(N)^{-1}$ establishes a graded vector space isomorphism:
\begin{equation}
    \bar{\Lambda}_N \simeq \left( V_N \cap \ker \partial_{p_1} \right)
\end{equation}
\end{lemma}
\begin{proof}
    On $V_N$ the operator $\hat{\sfC}(N)$ is invertible. At the same time, due to \eqref{eq:DualPieri} the operator $\partial_{p_1}$ preserves $V_N$. The equality then follows by direct evaluation on Schur polynomials with $\ell(\lambda)\leq N$:
    \begin{equation}
        \nabla_N \left( \dfrac{s_{\lambda}(\bx)}{\sfC_{\lambda}(N)} \right)= \sum_{\Box\in\Rem\lambda} \dfrac{s_{\lambda-\Box}(\bx)}{\sfC_{\lambda-\Box}(N)} = \pi_N \left(   \hat{\sfC}(N)^{-1} \sum_{\Box \in \Rem\lambda}s_{\lambda-\Box}\right)
    \end{equation}
\end{proof}
Let us summarize the idea again. We have found an explicit isomorphism between the space  $V_N \cap \mathbb{C}[p_2,p_3,\ldots]$ and the space of shift invariant symmetric polynomials $  \bar{\Lambda}_N$ by using the operator $\hat{\sfC}(N)$. Given a symmetric function that does not depend on the first power sum and has a length-bounded Schur expansion we can construct a shift invariant symmetric polynomial.

Let us reformulate this statement in another way that is useful later. Let $q \in \Lambda$ be a homogeneous degree $d$ symmetric polynomial and $q(\bx)=\pi_N q$, and denote its Schur expansion coefficients by $[q]_\lambda$. We have:
\begin{equation}
    [q]_\lambda=(q,s_\lambda)
\end{equation}
and 
\begin{equation}
    q(\bx)=\sum_{\substack{\lambda \vdash d\\\ell(\lambda)\leq N}} [q]_\lambda s_{\lambda}(\bx)
\end{equation}
The shift invariance condition is:
\begin{equation}
    \nabla_N q(\bx)=0
\end{equation}
Due to the lemma above this means that:
\begin{equation}
    \sum_{\Box \in \Add{\mu}}\sfC_{\mu+\Box}(N)[q]_{\mu+\Box} = 0 \,\quad \forall\mu \vdash d-1 \,. 
\end{equation}
Here the set of addable boxes is defined as:
\begin{equation}
    \Add\lambda = \left\{ (i,j) \in \mathbb{Z}_{>0}^2 \, : \, j=\lambda_{i}+1 \leq \lambda_{i-1}  \right\}
\end{equation} 
with, formally, $\lambda_{0}=\infty$. In other words, the set $\Add\lambda$ contains boxes $\Box$ at coordinates $(i_{\Box},j_{\Box})$ that can be added to the initial Young diagram $\lam$ in a way that the resulting diagram $\lam + \Box$ is a Young diagram. When $d\leq N$ we can also define the normalized coefficients:
\begin{equation}
    [q]'_{\lambda}= \sfC_\lambda(N)[q]_{\lambda}
\end{equation}
which should then satisfy the following recursion:
\begin{equation}\label{eq:recursion}
    \sum_{\Box \in \Add{\mu}}[q]'_{\mu+\Box} = 0 \,\quad \forall\mu \vdash d-1\,. 
\end{equation}
This equation looks very similar to relations satisfied by characters and related combinatorial quantities.

\section{The fusion basis}
\label{sec::forgotten}
Due to Lemma \ref{lemma:transfer} our goal is to find a basis in $V_N \cap \ker \partial_{p_1}$. We would like elements of this basis to be labeled by nuclear partitions in some natural way. Namely we will look for functions $q_{\lambda}$ such that  the condition $\partial_{p_1}q_{\lambda}=0$ would be solved by taking $\lambda$ to be nuclear. At the same time, we would like the functions to have Schur expansion coefficients that are either explicitly known or satisfy a recursion or can be calculated in some controlled way.

We find that these qualities are satisfied by the \textit{forgotten symmetric functions $f_{\lambda}$} \cite{Macdonald:book,Alexandersson2020}.
They are defined as the image of monomial symmetric functions under the involution $\omega$ \cite{Macdonald:book}:
\begin{equation}
f_{\lambda}=\omega(m_{\lambda})
\end{equation}
The involution acts on power sums as:
\begin{equation}
    \omega(p_k)=(-1)^{k-1}p_k\,.
\end{equation}
The Schur expansion follows from that of the monomial symmetric functions. 
\begin{proposition}
    The Schur expansion of the forgotten symmetric functions is given by:
     \begin{equation}
        f_\lambda=\sum_{ \mu^T\unlhd \lambda}(K^{-1})_{\lambda\mu^{T}} \, s_{\mu}
    \end{equation} 
    where $(K^{-1})$ is the inverse Kostka matrix. 
\end{proposition}
\begin{proof}
The Kostka matrix is defined as the transition matrix between the Schur basis and the monomial symmetric functions \cite{Macdonald:book}:
\begin{equation}
    s_{\lambda}=\sum_{\mu} K_{\lambda\mu} m_\mu
\end{equation}
The inverse transformation is therefore:
    \begin{equation}
        m_\lambda=\sum_{\mu \unlhd\lambda }K^{-1}_{\lambda\mu}s_\mu
    \end{equation}
 Now apply the involution to both sides. On the left we get the forgotten symmetric function. On Schur functions the involution $\omega$ acts by transposing the partition \cite{Macdonald:book}:
    \begin{equation}
        \omega(s_{\mu})=s_{\mu^T}\end{equation}
\end{proof}
The important corollary of this is that the lengths of the partitions indexing Schur polynomials in its expansion are controlled by the maximal component of $\lambda$:
\begin{corollary}\label{corollary:fbasis}
    For $\lambda$ such that $\lambda_i \leq N\,, \forall i$ we have $f_{\lambda} \in V_N$.
\end{corollary}
\begin{proof}
    This is implied by $\mu^T \unlhd \lambda \Rightarrow \mu \unrhd \lambda^T \Rightarrow \ell(\mu) \leq \ell(\lambda^T)=\lambda_1 \leq N$. Therefore in the Schur expansion only Schur polynomials with $\ell(\mu)\leq N$ contribute.
\end{proof}
The $p_1$ derivative acts on forgotten symmetric functions in a simple way.
\begin{proposition}
For forgotten symmetric functions one has:
     \begin{equation}
        \frac{\partial}{\partial p_1} f_{\lambda} = \left\{ \begin{array}{cc}
             f_{\lambda\backslash   (1)}\,, & m_1(\lambda) >0\\
           0\,,& m_1(\lambda)  = 0
        \end{array}\right.
    \end{equation}
where we used a common notation $\lambda = (1^{m_1(\lambda)}, 2^{m_2(\lambda)}, \ldots)$.
\end{proposition}
\begin{proof}
    By virtue of the involution one has $(e_{\mu},f_{\lambda})=\delta_{\lambda\mu}$. Then $(e_{\nu},\dfrac{\partial}{\partial p_1}f_{\lambda})=(p_1e_{\nu},f_{\lambda})=(e_{\nu\cup (1)},f_{\lambda})=\delta_{\lambda,\nu \cup (1)}$
\end{proof}
Combining this with Corollary \ref{corollary:fbasis} we have:
\begin{corollary}
    At fixed $N$ the forgotten symmetric polynomials with $1 <\lambda_i \leq N\,,\forall i$ form a basis for the subspace of $p_1$ independent polynomials in $V_N$:
    \begin{equation}
      V_N \cap  \ker  \partial_{p_1} = \operatorname{span}\left\{ f_{\lambda} \,: \lambda \in \mathcal{N}_N\right\}
    \end{equation}
\end{corollary}
Therefore, by Lemma \ref{lemma:transfer} we can construct a basis for the space of shift invariant polynomials, thereby obtaining the following theorem.
\begin{theorem}
    For $d=0$, set $W^{\rm f}_{\varnothing}(\bx)=1$, and for $d>0$, set
    $W^{\rm f}_{\lambda}(\bx)=\pi_N \hat{\sfC}(N)^{-1} f_{\lambda}$. Then for $\lambda \in \mathcal{N}_N$ the family of symmetric polynomials $W^{\rm f}_{\lambda}(\bx)$   is a  basis for the ring $\bar{\Lambda}_N$ and has the Schur expansion:
     \begin{equation}
        W^{\rm f}_{\lambda}(\bx)= \sum_{\substack{\mu^T\unlhd \lambda} } (K^{-1})_{\lambda\mu^T} \, \dfrac{ s_{\mu}(\bx)}{\sfC_{\mu}(N)}
    \end{equation}
\end{theorem}

\subsection{Uniformization in $N$}
\label{sec::uniform N}
Up to this stage we considered a basis in the space for each $N$ separately. The elements $W_{\lambda}^{\rm f}(\bx)$ are not stable under reducing the number of variables. Indeed, even the dimensions of the homogeneous components $\bar{\Lambda}_{N,d}$ change as we change $N$. For $N=4$ the partition $\lambda=[4]$ labels an independent basis element $W_{[4]}^{\rm f}(x_1,\ldots,x_4)$. Setting $x_4=0$ reduces the number of variables to $N=3$; however, $  W_{[4]}^{\rm f}(x_1,\ldots,x_3,0) \neq 0$. While, according to our definition, it is not present in the basis anymore as $\lambda_1>3$, the reduction does not happen naturally. This issue prevents identifying this basis in the stable ring.

In this section we fix this issue by defining elements in the stable ring that have the desired reduction. We choose a normalization of the basis elements such that the excess functions vanish automatically after the projection.

We start by defining the analog of the forgotten-like elements $W_{\lambda}^{\rm f}(\bx)$ in the stable ring of symmetric functions. Define the following elements $ \hat{\sfC}(u)^{-1} \cdot f_{\lambda}$ in the stable ring $\Lambda(u)=\Lambda \otimes \mathbb{C}(u)$. For a fixed integer $N$ such elements have the property that, if $\lambda_1\leq N$, then after projection $\pi_N(  \hat{\sfC}(u)^{-1} \cdot f_{\lambda})$ is finite at $u=N$ and $\hat{\sfC}(N)^{-1} \cdot f_{\lambda}$ is also finite before projection. However, for partitions with  $\lambda_1>N$ the evaluation $u=N$ and the projection do not commute. Namely, in this case, evaluation before projection is ill-defined. However, if we were to work polynomially in $u$, $\Lambda[u]=\Lambda \otimes \mathbb{C}[u]$, then this issue would not appear.
\begin{definition}
    Define $\tilde{\pi}_N$ as the map:
\begin{equation}
    \tilde{\pi}_N: \Lambda[u] \rightarrow \Lambda_N
\end{equation}
\begin{equation}
    q(u) \rightarrow \pi_N(q(N))
\end{equation}
\end{definition}
or equivalently 
\begin{equation}
    \tilde{\pi}_N q(u) = \pi_N q(u) \Big|_{u = N}
\end{equation}
In order to get rid of the denominators appearing in $W_{\lambda}^{\rm f}[u]$ we define their common multiple as:
\begin{equation}
    D_{\lambda}(u)=\underset{(K^{-1})_{\lambda\mu^T}\neq0}{ \operatorname{lcm}}  \sfC_{\mu}(u)
\end{equation}

Multiplying by these factors removes all denominators in $W_{\lambda}^{\rm f}[u]$.
    Define:
    \begin{definition}
            \begin{equation}
        \widetilde{W}^{\rm f}_{\lambda}[u]=D_{\lambda}(u)\hat{\sfC}(u)^{-1} \cdot f_{\lambda}
    \end{equation}
    \end{definition}
\begin{proposition}\label{prop::projection of fusion}
The projection $\tilde{\pi}_N$ of $\widetilde{W}^{\rm f}_{\lambda}[u]$ is well defined and is given by:
        \begin{equation}
           \tilde{\pi}_N(  \widetilde{W}^{\rm f}_{\lambda}[u]) = \left\{ \begin{array}{cc}
              D_{\lambda}(N) W^{\rm f}_{\lambda}(\bx)  & N \geq \lambda_1  \\
               0 & N < \lambda_1
           \end{array} \right.
        \end{equation}
\end{proposition}
\begin{proof}
    First, because all denominators have been cleared, $\widetilde{W}^{\rm f}_{\lambda}[u]$ is polynomial in $u$. Hence the projection and evaluation commute. 
    \\

    Let $N<\lambda_1$. The partitions $\mu$ in the Schur expansion satisfy $\ell(\mu) \leq \lambda_1$. Split them into two groups $\ell(\mu)\leq N$ and $\ell(\mu) > N$. Then for short partitions we have:
    \begin{equation}
       \dfrac{ {D}_{\lambda}(u)}{\sfC_{\mu}(u)}\Bigg|_{u=N}=0 \,,\, \text{if } \ell(\mu) \leq N 
    \end{equation}
    since $\sfC_{\mu}$ cancels in $D_{\lambda}$, but it does not cancel $\sfC_{\lambda^{\rm T}} | D_{\lambda}$ and $\sfC_{\lambda^{\rm T}} (N)=0$.
    Therefore the corresponding terms in the Schur expansion vanish.
For longer partitions the coefficient is polynomial and hence:
    \begin{equation}
         \dfrac{ {D}_{\lambda}(u)}{\sfC_{\mu}(u)}\Bigg|_{u=N} <\infty \,,\, \text{if } \ell(\mu) > N\,.
    \end{equation}
    However, $\pi_N(s_{\mu})=0 $ for $\ell(\mu) > N$. Therefore these terms also vanish because of the projection.
    \\
    When $\lambda_1\leq N$ all denominators $\sfC_{\mu}(N)0$ in the Schur expansion are nonzero since $\ell(\mu)\leq N$ and $D_{\lambda}(N)\neq 0$.
\end{proof}

Therefore the elements $\widetilde{W}^{\rm f}_{\lambda}[u]$ in the ring of symmetric functions interpolate between shift invariant bases at different $N$. Define $ \bar{\Lambda}[u]$ to be a subset of the stable ring $\Lambda[u]$ that is spanned by functions $\widetilde{W}^{\rm f}_{\lambda}[u]$ for nuclear partitions:
\begin{equation}
   \bar{\Lambda}[u]= \operatorname{span}_{\mathbb{C}[u]}\left\{ \widetilde{W}^{\rm f}_{\lambda}[u]   \ : \ \lambda \in \mathcal{N}\right\}
\end{equation}
\begin{proposition}
The projection $\tilde{\pi}_N$ is a surjective vector space homomorphism:
    \begin{equation}
        \bar{\Lambda}[u] \rightarrow \bar{\Lambda}_N
    \end{equation}
\end{proposition}
This space should play a role similar to that of the stable ring for shift invariant symmetric functions. A more precise qualification is beyond the scope of this paper.

\subsection{Examples}\label{sec::fusion examples}
Let us provide a few examples. At degrees two and three there is only one basis element. Therefore we start at degree 4. For $\lambda=[4]$ we get:
\begin{equation}
\begin{aligned}
    \widetilde{W}^{\rm f}_{[4]}[u] = -(u-3) (u-2) (u-1) \, s_{[4]}+(u-3) (u-2) (u+3) \, s_{[3,1]} \\
    -(u-3) (u+2) (u+3) \, s_{[2,1,1]}+(u+1) (u+2) (u+3) \, s_{[1,1,1,1]}
\end{aligned}
\end{equation}
and for $\lambda=[2,2]$ we get:
\begin{equation}
    \widetilde{W}^{\rm f}_{[2,2]}[u] = (u-1) u \, s_{[4]} - (u+3) u \, s_{[3,1]} + (u+2) (u+3) \, s_{[2,2]}
\end{equation}
Let us explain how these two elements behave under projection. Both polynomials are non-zero under $\tilde{\pi}_{k}$, when $k \geq 4$. This is in agreement with the fact that at $N \geq 4$ there should be two linearly independent shift invariant polynomials. At $\tilde{\pi}_3$ the first polynomial $ \widetilde{W}^{\rm f}_{[4]}[u] $ vanishes. At $u = 3$, the first three terms vanish due to vanishing coefficients and the last term vanishes because $\pi_3 s_{[1,1,1,1]} = 0$. Similarly, this polynomial is zero under the projector $\tilde{\pi}_{2}$. At the same time the second element $\widetilde{W}^{\rm f}_{[2,2]}[u]$ is non-vanishing, both under $\tilde{\pi}_{2}$ and $\tilde{\pi}_{3}$. It gives the degree four shift invariant polynomial at two and three variables, respectively.
\subsection{Properties}
We collect some of the properties of the ingredients of the construction: the Kostka coefficients and the forgotten symmetric functions.
The elements of the basis are defined explicitly in terms of the inverse Kostka matrix. The main advantage of such a representation is that Kostka coefficients have been intensively studied and therefore many of their properties are known. In particular both the direct and inverse Kostka coefficients can be computed via recursive formulas. Here we provide one of them as an illustration. 
\begin{theorem}[Eğecioğlu–Remmel \cite{e1990combinatorial}]
Let
\[
\lambda=(1^{m_1},2^{m_2},\ldots,d^{m_d}),
\qquad
\mu=(\mu_1,\ldots,\mu_\ell)\vdash d,
\qquad
\mu_1\ge\cdots\ge\mu_\ell>0.
\]
For \(r\ge1\) with \(m_r>0\), let \(\lambda\setminus r\) denote the partition obtained by removing one part of size \(r\); equivalently, its multiplicities are \(m_j-\delta_{jr}\).

For \(1\le i\le\ell\), define
\[
r_i=\mu_i+\ell-i,
\qquad
\mu^{(i)}
=
(\mu_1,\ldots,\mu_{i-1},
 \mu_{i+1}-1,\ldots,\mu_\ell-1),
\]
with zero parts omitted. Then
\[
(K^{-1})_{\lambda\mu}
=
\sum_{\substack{1\le i\le\ell\\m_{r_i}>0}}
(-1)^{\ell-i}
(K^{-1})_{\lambda\setminus r_i,\,\mu^{(i)}}.
\]
The initial condition is
\[
(K^{-1})_{\varnothing,\varnothing}=1,
\]
and entries indexed by partitions of unequal sizes are understood to be zero.
\end{theorem} 
The recursion expresses the inverse Kostka coefficient for $\lambda,\mu$ in terms of coefficients for smaller partitions. For example, at degree four, we have:
\begin{equation}
    (K^{-1})_{[4],[3,1]}= - (K^{-1})_{\varnothing,\varnothing}
\end{equation}
since $r_1=4,r_2=1$. For the partition $\lambda=[2,2]$ we get:
\begin{equation}
  \begin{gathered}     (K^{-1})_{[2,2],[2,2]}=  (K^{-1})_{[2],[2]}\,, \quad   (K^{-1})_{[2,2],[2,1,1]}=-  (K^{-1})_{[2],[2]}\,, 
        \\
        (K^{-1})_{[2,2],[1^4]}=- (K^{-1})_{[2],[1,1]}\,.
  \end{gathered}
\end{equation}
There are several useful formulas that allow one to compute forgotten symmetric functions. Though they are counterparts of the respective formulas for the monomial functions, we list them here. While their Schur expansion is given by inverse Kostka coefficients, which can be computed recursively, their power sum expansion is governed by more intuitive combinatorics. First, the expansion can be packed into a generating function; see Appendix \ref{sec::append A} for details. Let $t_{k} \, ,k\geq 1$ be auxiliary power sum-like variables. Then:
    \begin{equation}
    \label{gen funct for forgotten}
        \exp\left( \sum_{\substack{\lambda \in \mathcal{P}
        \\
        \lambda \neq \varnothing}} (-1)^{|\lam| + \ell(\lam) }v_{\lambda} \, p_{|\lambda|} \cdot t_{\lambda}  \right) = \sum_{\mu \in \mathcal{P}} f_{\mu} t_{\mu}
    \end{equation}
    where $t_{\lambda} = \prod_i t_{\lambda_i}$ and  
    \begin{equation}
        v_{\lambda} = \frac{(\ell(\lam) - 1)!}{a_{\lambda}}
    \end{equation}
where $a_{\lambda}=\prod\limits_{j \geq 1} m_j(\lam)!$.
By inspecting this expansion one gets the following formula, which motivates the term \emph{fusion}:
    \begin{equation}
        f_{\mu} = \epsilon_\mu\sum_{s = 1}^{\infty} \frac{1 }{s!}\sum_{\substack{ \lambda^{(1)} \cup \lambda^{(2)} \ldots \cup \lambda^{(s)} : = \mu \\ \ell(\lam^{(i)}) \geq 1 }}  \prod\limits_{i = 1}^{s} v_{\lam^{(i)}} p_{|\lam^{(i)}|}
    \end{equation}
    where \[
\epsilon_\mu=(-1)^{|\mu|+\ell(\mu)}.
\]
Such a formula was originally established in \cite{Merca2015}. Here we write it in slightly different notation. The sum is over all possible ways to divide the partition $\mu$ into non-overlapping subsets. For each subset $\lambda^{(i)}$ the respective power sum $p_{|\lam^{(i)}|}$ depends only on the total sum of the elements of the set. Let us explain how this can be understood as ``fusion'', since this is how we originally came up with this construction before we were able to identify the relevant basis as the forgotten symmetric functions. For a fixed $\mu$, take its division into a maximal number of subsets, that is, each subset is simply a single component of $\mu$ and $s=\ell(\mu)$. There are precisely $\frac{\ell(\mu)!}{a_\mu}$ such divisions, which give the same term in the sum -- $p_\mu$. We can simply choose the representative:
\begin{equation}
    \lambda^{(i)}=[\mu_i]\,.
\end{equation} 
Therefore $f_\mu = \dfrac{\epsilon_\mu}{a_\mu} p_\mu+(\text{terms with  $s<\ell(\mu)$}) $. Next, the term with $s=\ell(\mu)-1$ will contain divisions with a single subset containing two elements $\mu_{k},\mu_l$. The corresponding contribution looks like:
\begin{equation}
    p_{\mu_{k}+\mu_l}\prod_{\substack{i=1 \\ i \neq k,l}}^{\ell(\mu)}p_{\mu_i}
\end{equation}
The combinatorial coefficients  will then count the number of ways to obtain such terms. It looks like $\mu_k$ and $\mu_l$ have fused. Take an example:
\begin{equation}
    f_{[3,2]}=-p_{2}p_3-p_5
\end{equation}
Here the initial $3$ and $2$ components have fused to produce $5$. For longer partitions:
\begin{equation}
          f_{[2,1,1,1]} = -\frac{1}{6} \left(p_1^3 p_2+3 p_1 p_2^2+3 p_1^2 p_3+5 p_2 p_3+6 p_1 p_4+6 p_5\right)
\end{equation}
The first few fusions are:
\begin{equation}
   \begin{array}{ccccc}
    [2,1,1,1] &\rightarrow &[2,1+1,1] &=&[2,2,1]
    \\
        {[2,1,1,1]} &\rightarrow &[2+1,1,1] &=&[3,1,1]
        \\
        {[2,1,1,1]} &\rightarrow & \left\{\begin{gathered}
             {[2+1,1+1]}\\
              {[2,1+1+1]}
        \end{gathered}\right\}
         &=& [3,2]
    \end{array}
\end{equation}
\section{The character, decay and fission  bases}

Many of the complications encountered in the previous section were due to the kernel of the $\hat{\sfC}$ operator and the projection $\pi_N$. For fixed $N$ and $d\leq N$ the kernel is trivial. Moreover in this range the space of shift invariant symmetric polynomials is labeled by nuclear partitions without any restrictions on the maximal component. In this section we first study that the subspace $\bar{\Lambda}_{N,\leq d}$, which has other natural bases. For the first two subsections, unless otherwise specified, we assume that $d\leq N$ and $|\lambda|\leq N$. Next, we use the results obtained for the decay construction to present another basis, which we call the fission basis, and prove its shift invariance. The fission basis, akin to the fusion basis, is well defined in any degree $d$.

\subsection{Character basis}
Let us work within the same construction as before and rely on Lemma \ref{lemma:transfer}. For degree $d\leq N$ we are always in $V_N$: the span of Schur polynomials for partitions with lengths $\ell(\mu)\leq N$. 
The most natural functions independent of $p_1$ that span $\ker \partial_{p_1}$ are simply the power sum monomials  $p_\mu$ for nuclear partitions $\mu \in \mathcal{N}$. Their expansion in terms of Schur functions is straightforward and follows from standard Frobenius formulas, namely:
\begin{equation}
    p_\mu = \sum_{\lambda \vdash |\mu|} \chi_{\mu}^{\lambda} \, s_{\lambda}
\end{equation}
where $\chi_{\mu}^{\lambda}$ are symmetric group characters. The requirement that $\mu$ be nuclear means that  $\chi_{\mu}^{\lambda}$ are restricted to conjugacy classes of derangements - permutations that do not contain length one cycles. Hence
\begin{equation}
      W_{\mu}^{\rm p}(\bx)= \pi_N \hat{\sfC}^{-1}(N)  p_\mu \,\quad \mu \in \mathcal{N}
\end{equation}
which can be explicitly written in terms of Schur polynomials as:
\begin{equation}
    W_{\mu}^{\rm p}(\bx)=\sum_{\lambda \vdash |\mu|} \chi_{\mu}^{\lambda} \, \dfrac{s_{\lambda}(\bx)}{\sfC_{\lambda}(N)} 
\end{equation}
Since $d\leq N$ the finite-variable projection is well defined, and we have the following result.
\begin{proposition}
    Symmetric polynomials  $W_{\mu}^{\rm p}(\bx)$  for $\mu \in \mathcal{N}(d)$ are a vector space basis for the space of shift invariant symmetric polynomials of degree $d \leq N$.
\end{proposition}
We slightly abuse the notation and use the same symbol for the element of $\Lambda$ and its image in $\Lambda_N$.

The main drawback of this basis is that its expansion in terms of Schur functions is not triangular. While this is not a property that we required initially, it is rather useful in practice. Namely, given a combination of Schur polynomials that is known to be shift invariant, one cannot easily read off its expansion in terms of the character basis without referring to the character table.

\subsection{The decay basis}
\label{sec::decay basis}
In this section we construct the decay basis, which was first proposed in \cite{Dixon:2026fuz}. Similar constructions also appeared in \cite{deMelloKoch:2004crq, Mishnyakov:2020khb}. It is still only defined when $d \leq N$, however there are two important distinctions from the character basis. First, as we will see, it has a triangular expansion in terms of Schur functions. Second, the coefficients of its Schur expansion can be computed in closed form by Aitken's determinant formula, as opposed to symmetric group characters.

Let us first review how it was originally constructed \cite{Dixon:2026fuz}. For a partition $\nu$ of length $l$ use the notation $\nu  \cup (1^k) $ for a partition of $|\nu|+k$ with components $[\nu_1,\ldots,\nu_l,\underbrace{1,\ldots, 1}_{k}]$. Then, for each nuclear partition $\lambda$, set:
\begin{equation}
    W_{\lambda}^{\rm d}(\bx)= \sum_{\substack{\nu \subseteq\lambda\\
    \ell(\nu)=\ell(\lambda)}} (-1)^{|\lambda|-|\nu|}f^{\lambda/\nu} \dfrac{s_{\nu  \cup (1^{|\lambda|-|\nu|})}(\bx)}{\sfC_{\nu  \cup (1^{|\lambda|-|\nu|})}(N)}
\end{equation}

\paragraph{Examples:}
At degree $2$ this construction gives the same polynomial as the fusion basis up to a sign. The first non-trivial examples occur at $d=4$ where there exist two basis polynomials, corresponding to nuclear partitions $[4]$ and $[2,2]$. For $\lambda=[4]$, the relevant partitions $\nu$ are $[4],[3],[2],[1]$. Therefore for the full invariant we have:
\begin{equation}
 \begin{split}
        W_{[4]}^{\rm d}(\bx)= &\dfrac{s_{[4]}(\bx)}{N(N+1)(N+2)(N+3)}-\dfrac{s_{[3,1]}(\bx)}{N(N+1)(N+2)(N-1)}\\+&\dfrac{s_{[2,1,1]}(\bx)}{N(N+1)(N-1)(N-2)}-\dfrac{s_{[1,1,1,1]}(\bx)}{N(N-1)(N-2)(N-3)}
 \end{split}
\end{equation}
For $\lambda=[2,2]$, the subpartitions $\nu$ are $[2,2],[2,1]$ and $[1,1]$, which gives the expansion:
\begin{equation}
 \begin{split}
        W_{[2,2]}^{\rm d}(\bx)= &\dfrac{s_{[2,2]}(\bx)}{N^2(N+1)(N-1)}\\-&\dfrac{s_{[2,1,1]}(\bx)}{N(N+1)(N-1)(N-2)}+\dfrac{s_{[1,1,1,1]}(\bx)}{N(N-1)(N-2)(N-3)}
 \end{split}
\end{equation}
 \begin{lemma}[\cite{Dixon:2026fuz}]
    \label{lemma::decay shift invariance proof}
    Symmetric polynomials $ W_{\lambda}^{\rm d}(\bx)$ for $\lambda \in \mathcal{N}$ form a vector space basis of $\bar{\Lambda}_{N, \leq N}$
\end{lemma}
\begin{proof}
    The linear independence of these polynomials follows from the Schur expansion. Indeed, in each sum, the only Schur polynomial corresponding to a nuclear partition is $s_{\lambda}(\bx)$. Now let us see why they are shift invariant.  For this consider equation \eqref{eq:recursion}:
    \begin{equation}
        \left[W_{\lambda}^{\rm d}(\bx) \right]'_\mu= \left\{ \begin{array}{cc}
            (-1)^{k}f^{\lambda/\nu} \,, & \mu=\nu \cup(1^{k}) \ \text{for}\ \ell(\nu)=\ell(\lambda) \\
            0\,, &  \text{otherwise}
        \end{array} \right.
    \end{equation}
Let $\mu$ be a partition of $d-1$. Then in \eqref{eq:recursion} the only non-zero terms are such that:
\begin{equation}
    \mu+\Box = \nu(\Box) \cup(1^{k(\Box)})
\end{equation}
By going over the possible cases one can see that there are two types of non-zero contributions. One, with $\nu(\Box)=\nu+\Box$ where the box is added to one of the nonzero components of $\nu$ and $k(\Box)=k$ and the other with $\nu(\Box)=\nu$ and $k(\Box)=k+1$, where the box is added to the column of $1$'s. Then \eqref{eq:recursion} turns into:
\begin{equation}
    f^{\lambda/\nu}-\sum_{i=1}^{\ell(\lambda)}
f^{\lambda/(\nu+\Box_i)}=0\,,
\end{equation}
which is satisfied identically.
\end{proof}
The special feature of this basis is that the coefficients $f^{\lambda/\nu}$ have a simple closed form expression, in contrast to Kostka numbers or characters of symmetric groups which are computationally more involved. It is known as Aitken's determinantal formula (Corollary 7.16.3 in ref.~\cite{stanley2023enumerative2}):
\begin{equation}\label{eq:flmdetapp}
        f^{\lambda/\mu} = (|\lambda|-|\mu|)!  \,\det \left( \dfrac{1}{(\lambda_i-\mu_j+j-i)!} \right) \,,
    \end{equation}
where $1/k!$ is defined to be $0$ if $k<0$ for $k=\lambda_i-\mu_j+j-i$.
\\

It is interesting to study the preimages of the decay basis elements under the $\pi_N \hat{\sfC}(N)^{-1}$ map. Let, for any $|\lambda|=d$, without the restriction:
\begin{equation}
r_\lambda=
\sum_{\substack{\nu\subseteq\lambda\\
                 \ell(\nu)=\ell(\lambda)}}
(-1)^{d-|\nu|}f^{\lambda/\nu}
s_{\nu\cup(1^{d-|\nu|})},
\end{equation}
When $d\leq N$, this symmetric function agrees with the preimage of the decay basis, namely:
\begin{equation}
  \pi_N \hat{\sfC}(N)^{-1}   r_{\lambda}= W_{\lambda}^{\rm d}\,.
\end{equation}
It turns out that they have a special expansion in terms of power sums.
\begin{lemma}
    Let $\Xi_{\lambda}^{\mu}$ denote the entries of the inverse of the restriction of the character table to nuclear partitions.
    Then for $\lambda \in \mathcal{N}(d)$:
    \begin{equation}
          r_{\lambda}=\sum_{\mu \in \mathcal{N}(d)} \Xi_{\lambda}^\mu \, p_\mu
    \end{equation}
    and the inverse transform:
    \begin{equation}
    p_{\lambda}=\sum_{\mu \in \mathcal{N}(d)} \chi_{\lambda}^\mu \, r_\mu
    \end{equation}
\end{lemma}
\begin{proof}
Let
\begin{equation}\label{eq:rA}
    r_{\lambda}=\sum_{\mu \vdash d} A_{\lambda}^{\mu} \, p_\mu
\end{equation}
    The proof of shift invariance of $W_{\lambda}^d$ (see Lemma \ref{lemma::decay shift invariance proof}) can be extended straightforwardly to give:
   \begin{equation}
      \dfrac{\partial}{\partial p_1} r_{\lambda}=0
   \end{equation}
   Hence the sum in \eqref{eq:rA} is actually restricted to nuclear partitions:
   \begin{equation}
    r_{\lambda}=\sum_{\mu \in  \mathcal{N}(d)} A_{\lambda}^{\mu} \, p_\mu
\end{equation}
    By definition of $r_{\lambda}$ the only Schur function corresponding to a nuclear partition in its expansion is $s_{\lambda}$. Hence, for nuclear partitions $r_{\lambda}$ are orthogonal to Schur functions with respect to to the Hall inner product:
    \begin{equation}
(r_{\lambda},s_{\nu})=\delta_{\lambda\nu} \,,\quad \lambda,\nu  \in \mathcal{N}(d)
    \end{equation}
    We can find $\Xi$ by simply rewriting this condition in the power sum basis:
    \begin{equation}
         \sum_{\mu\in \mathcal{N}(d)} A_{\lambda}^{\mu} (p_\mu,s_\nu)=\sum_{\mu\in \mathcal{N}(d)} A_{\lambda}^{\mu} \chi^{\nu}_\mu =\delta_{\lambda\nu} \, ,\quad \lambda,\nu  \in \mathcal{N}(d)
    \end{equation}
    By construction we get:
    \begin{equation}
A_{\lambda}^{\mu}=\Xi_{\lambda}^\mu\,.
    \end{equation}
\end{proof}

For example in degree four the restricted character matrix is given by:
\begin{equation}
    [\chi]\Big|_{\mathcal{N}(4)}= \left(\begin{array}{cc}
        1 &0  \\
         1& 2
    \end{array}\right)
\end{equation}
and hence:\\\begin{equation}
\Xi=\left(\begin{array}{cc}
        1 &0  \\
         -\frac12&\frac12
    \end{array}\right)
\end{equation}
This gives:
\begin{equation}
    r_{[4]}=p_4\,,\quad r_{[2,2]}=\dfrac{p_2^2-p_4 }{2}
\end{equation}
\begin{corollary}
    For $d\leq N$ we have transition formulas:
    \begin{equation}
        \Wp_{\mu}(\bx)=\sum_{\lambda \in \mathcal{N}(d)} \chi^{\lambda}_\mu \, \Wd_\lambda(\bx)
    \end{equation}
    and
    \begin{equation}
        \Wf_{\mu}(\bx)=\sum_{\lambda  \in \mathcal{N}(d)} (K^{-1})_{\mu\lambda^T} \, \Wd_\lambda(\bx)
    \end{equation}
\end{corollary}
\begin{proof}
    For the character basis the statement follows directly. For the fusion basis, recall
    \begin{equation}
        (s_\mu,r_\lambda)=\delta_{\mu\lambda} 
    \end{equation}
    for $\mu \in \mathcal{N}(d)$. Let:
    \begin{equation}
        f_\mu=\sum_{\nu \in \mathcal{N}(d)}\alpha_{\mu\nu}r_\nu
    \end{equation}
    By pairing this with Schur functions we get:
    \begin{equation}
        (s_{\lambda} ,f_\mu)=\sum_{\nu \in \mathcal{N}(d)}\alpha_{\mu\nu} (s_{\lambda},r_\nu)
    \end{equation}
    Hence
    \begin{equation}
        \alpha_{\mu\lambda}= (K^{-1})_{\mu \lambda^{\rm T}}\,.
    \end{equation}
\end{proof}

\subsection{The fission basis}
Based on the decay basis, we can construct another family of  shift invariant symmetric functions which can now be extended to every degree $d$. Hence, below we will lift the restriction from the beginning of this section. The construction is based on the property of the $p_1$ derivative with respect to the involution.

\begin{lemma}
    The subspace $\ker \partial_{p_1}$ is preserved by the involution $\omega$. In other words if $\frac{\partial}{\partial p_1} a = 0$ for $a \in \Lambda$ then
    \begin{equation}
        \frac{\partial}{\partial p_1} \omega(a)=0
    \end{equation}
\end{lemma}
\begin{proof}
    The involution is self dual with respect to the Hall inner product, hence the chain of equalities:
    \begin{equation}
        (\omega(b),\partial_{p_1} a)= (p_1\,\omega(b),a)=(\omega(p_1\,b),a)=(p_1 \,b,\omega(a))=(b,\partial_{p_1}\omega(a))
    \end{equation}
\end{proof}

We apply the involution to the preimages of the decay basis to obtain preimages of the new basis that we call fission
\begin{equation}
    g_{\lam} := \omega(r_{\lam}) = \sum_{\substack{\nu \subseteq\lambda\\
    l(\nu)=l(\lambda)}} (-1)^{|\lambda|-|\nu|}f^{\lambda/\nu} s_{(\nu  \cup (1^{|\lambda|-|\nu|}))^{T}}
\end{equation}
We can rewrite the formula for polynomials $g_{\lam}$ in terms of transposed partitions in order to obtain the form equivalent to \eqref{fission def}
\begin{equation}
        g_\lambda=\sum_{\substack{\eta \subseteq\lambda^{\rm T}\\
    \eta_1=\ell(\lambda)}} (-1)^{|\lambda|-|\eta|}f^{\lambda^{\rm T}/\eta} s_{\eta^{(|\lambda|)}}
    \end{equation}
    Indeed set $\eta=\nu^T$, then $\eta^{(k)}=\left(\nu^{[k]} \right)^{\rm T}$. In the formula above we use property $f^{\lam / \nu} = f^{\lam^{T} / \nu^{T}}$. According to this definition we have the following statement
\begin{proposition}
    At fixed $N$ the functions $g_{\lam}$ with $1 <\lambda_i \leq N\,,\forall i$ form a basis for the subspace of $p_1$ independent functions in $V_N$:
    \begin{equation}
      V_N \cap  \ker  \partial_{p_1} = \operatorname{span}\left\{ g_{\lambda} \,: \lambda \in \mathcal{N}_N\right\}
    \end{equation}
\end{proposition}
\begin{proof}
    There are precisely $|\mathcal{N}_N(d)|$ functions $g_\lambda$ in each degree. Their linear independence follows from the Schur expansion. Length of a Young diagram involved in the expansion is bounded 
    \begin{equation}
        l((\nu  \cup (1^{|\lambda|-|\nu|}))^{T}) = \nu_1 \leq \lam_1 \leq N
    \end{equation}
    That follows from $\nu \subseteq \lam$. Therefore we have $g_{\lam} \in V_N$. Together with $\frac{\partial r_{\lam}}{\partial p_1} = 0, \lam \in \mathcal{N}$ we have $\frac{\partial g_{\lam}}{\partial p_1} = 0$.
\end{proof}
By Lemma \ref{lemma:transfer} we construct a basis for the space of shift invariant polynomials
\begin{theorem}
    For $d=0$, set $W^{\rm F}_{\varnothing}(\bx)=1$, and for $d>0$, set
    $W^{\rm F}_{\lambda}(\bx)=\pi_N \hat{\sfC}(N)^{-1} g_{\lambda}$. Then for $\lambda \in \mathcal{N}_N$ the family of symmetric polynomials $W^{\rm F}_{\lambda}(\bx)$  is a  basis for the ring $\bar{\Lambda}_N$ and has the Schur expansion:
     \begin{equation}
        W^{\rm F}_{\lambda}(\bx)= \sum_{\substack{\nu\subseteq \lambda\\
        l(\nu)=l(\lambda)}} (-1)^{|\lambda|-|\nu|} \, f^{\lambda/\nu} \, \dfrac{s_{(\nu^{[|\lambda|]})^T}(\bx)}{\sfC_{(\nu^{[|\lambda|]})^{T}}(N)}
    \end{equation}
    where $\nu^{[|\lambda|]} = \nu  \cup (1^{|\lambda|-|\nu|})$.
\end{theorem}
At the next step we want to promote this basis to a basis in the stable ring and carry out a procedure of uniformization in $N$ as in sec.\ref{sec::uniform N}. We do not repeat the whole proof here but go through the key arguments. First we define the common multiple 
\begin{equation}
    D^{\prime}_{\lambda}(u)=\underset{\substack{\nu\subseteq \lambda\\ l(\nu)=l(\lambda) \\ f^{\lambda/\nu}\neq0} \\ }{ \operatorname{lcm}}  \sfC_{(\nu^{[|\lambda|]})^{T}}(u)
\end{equation}
then we introduce the following
\begin{definition}
            \begin{equation}
         \widetilde{W}^{\rm F}_{\lambda}[u]=D^{\prime}_{\lambda}(u)\hat{\sfC}(u)^{-1} \ g_{\lambda}
    \end{equation}
\end{definition}
\begin{proposition}
The projection $\tilde{\pi}_N$ of $\widetilde{W}^{\rm F}_{\lambda}[u]$ is well defined and is given by:
        \begin{equation}
           \tilde{\pi}_N(  \widetilde{W}^{\rm F}_{\lambda}[u]) = \left\{ \begin{array}{cc}
              D^{\prime}_{\lambda}(N) W^{\rm F}_{\lambda}(\bx)  & N \geq \lambda_1  \\
               0 & N < \lambda_1
           \end{array} \right.
        \end{equation}
\end{proposition}
We do not repeat the proof here as it follows that of Proposition \ref{prop::projection of fusion}.
\paragraph{Examples}: for the two nuclear partitions of degree four, $[4]$ and $[2,2]$, the fission polynomials coincide with the fusion polynomials provided in sec.\ref{sec::fusion examples}. The first difference occurs in degree six:
\begin{equation}
\begin{aligned}
\widetilde W^{\mathrm F}_{[3,3]}[u]
={}& 2u(u-2)(u-1)^2\,s_{[6]}
     -2u(u-2)(u-1)(u+5)\,s_{[5,1]}\\
 & +(u-2)(u-1)(u+4)(u+5)\,s_{[4,2]} +u(u-1)(u+4)(u+5)\,s_{[4,1,1]}\\
 & -(u-1)(u+3)(u+4)(u+5)\,s_{[3,2,1]} +(u+2)(u+3)(u+4)(u+5)\,s_{[2,2,2]} 
\end{aligned}
\end{equation}
\begin{equation}
\begin{aligned}
\widetilde W^{\mathrm F}_{[2,2,2]}[u]
={}&-u(u-1)(u+1)\,s_{[6]}
     +u(u+1)(u+5)\,s_{[5,1]}\\
 & -(u+1)(u+4)(u+5)\,s_{[4,2]}
     +(u+3)(u+4)(u+5)\,s_{[3,3]}
\end{aligned}
\end{equation}
The comparison with fusion in this two-dimensional span takes the form
\begin{equation}
\begin{aligned}
 \widetilde W^{\mathrm f}_{[2,2,2]}[u]
 &=\widetilde W^{\mathrm F}_{[2,2,2]}[u],\\
 \widetilde W^{\mathrm f}_{[3,3]}[u]
 &=(u+1)\widetilde W^{\mathrm F}_{[3,3]}[u]
   +(u-2)(u-1)\widetilde W^{\mathrm F}_{[2,2,2]}[u]
\end{aligned}
\end{equation}
Indeed, the corresponding preimages satisfy
$f_{[2,2,2]}=g_{[2,2,2]}$ and $f_{[3,3]}=g_{[3,3]}+g_{[2,2,2]}$.
Under $\widetilde\pi_{\geq 3}$ both elements survive. Under $\widetilde\pi_2$ the second element vanishes. Indeed, each coefficient of a diagram with at most two rows contains the factor $(u-2)$, and all
remaining Schur functions vanish in two variables. For the second element we obtain
\begin{equation}
 \widetilde\pi_2\bigl(\widetilde W^{\mathrm F}_{[2,2,2]}[u]\bigr) =-6(x_1-x_2)^6 
\end{equation}

\section{Applications}

\subsection{Specialization of the Harish-Chandra–Itzykson–Zuber integral}

Consider the following unitary random matrix integral \cite{harish1957differential,itzykson1980planar}:
\
\begin{equation}
    \sfZ_N^{\rm IZ}(X,Y) = \dfrac{1}{\operatorname{Vol}(U(N))} \int\limits_{U(N)}  [dU] \exp\left( \Tr UX U^{\dagger}Y \right)
\end{equation}
where $X$ and $Y$ are Hermitian matrices. It is well known that the integral can be evaluated and expressed in terms of eigenvalues of $X$ and $Y$:
\begin{equation}\label{eq:IZresult}
     \sfZ_N^{\rm IZ}(X,Y) = \left(\prod_{j=1}^{N-1}j! \right)\dfrac{\det e^{x_i y_j} }{\Delta(\by)\Delta(\bx)}
\end{equation}
Since the function is symmetric in both $x_i$ and $y_i$ one can expect a series expansion in Schur polynomials. It is well known, based either on the original integral or \eqref{eq:IZresult} that the result can be expanded into a convergent series \cite{Balantekin:2000vn}:
\begin{equation}\label{eq:IZchar}
      \sfZ_N^{\rm IZ}(X,Y) = \sum_{\lambda:\, \ell(\lambda)\leq N} \dfrac{s_{\lambda}(\bx) s_{\lambda}(\by)}{\sfC_{\lambda}(N)}
\end{equation}
This clearly resembles the expansion of shift invariant polynomials in terms of Schur functions. The connection can be made as follows. 
Under a diagonal shift of the $X$ variable the partition function transforms as:
\begin{equation}
    \sfZ_N^{\rm IZ}(X+a I,Y) = e^{a \Tr Y} \sfZ_N^{\rm IZ}(X,Y)
\end{equation}
Therefore for $\Tr Y=0$ the whole function is shift invariant. Another way to understand this property is to recall that the HCIZ integral is known to solve the bi-spectral problem of the Calogero-Moser system at the free fermion point \cite{calogero1969ground,bergere2009some,Morozov:2009jv}. Denote the Hamiltonians of this system as $\mathcal{H}_k^{\rm CM}(x)$. In particular we have:
\begin{equation}
    \mathcal{H}_1^{\rm CM}(x) = \sum_{i=1}^N \dfrac{\partial}{\partial x_i}
\end{equation}
and
\begin{equation}
     \mathcal{H}_2^{\rm CM}(x) = \sum_{i=1}^N \dfrac{\partial^2}{\partial x_i^2} +2 \sum_{i<j} \dfrac{1}{(x_i - x_j)} \left( \dfrac{\partial}{\partial x_i} - \dfrac{\partial}{\partial x_j}\right)
\end{equation}
With respect to the Calogero Hamiltonian it exhibits the property of bi-spectrality:
\begin{equation}
  \begin{split}
        \mathcal{H}_k^{\rm CM}(\bx)\sfZ_N^{\rm IZ}(X,Y)=p_k(\by) \sfZ_N^{\rm IZ}(X,Y)
        \\
          \mathcal{H}_k^{\rm CM}(\by)\sfZ_N^{\rm IZ}(X,Y)=p_k(\bx) \sfZ_N^{\rm IZ}(X,Y)
  \end{split}
\end{equation}
Therefore setting $y_i$ to be such that $p_1(\by)=0$ we get that the corresponding HCIZ integral is shift invariant. Apart from this condition, this statement is true for any values of the $\by$ parameters. Hence we can treat the corresponding function as a generating function of shift invariant polynomials. Therefore by construction we have:
\begin{proposition}
    The HCIZ integral with $\Tr Y=0$ is a generating function for the fusion basis $W_{\lambda}^{\rm f}(\bx)$. 
\end{proposition}
\begin{proof}
    It follows from the above discussion by reexpanding \eqref{eq:IZchar} into  elementary symmetric functions of $\by$. The forgotten symmetric functions are orthogonal to the elementary symmetric polynomials:
    \begin{equation}(f_{\lambda},e_\mu)=\delta_{\lambda\mu}
    \end{equation}
Therefore, one can reorganize the series \eqref{eq:IZchar} as:
    \begin{equation}
       \begin{split}
            &\sfZ_N^{\rm IZ}(X,Y)= \hat{\sfC}_x^{-1}(N)\sum_{d=0}^{\infty}  \sum_{\substack{\lambda \vdash d\\\ell(\lambda)\leq N} } s_{\lambda}(\bx)s_{\lambda}(\by)=
            \\&\hat{\sfC}_x^{-1}(N)\sum_{d=0}^{\infty} \sum_{\substack{\lambda \vdash d\\\lambda_1 \leq N} } f_{\lambda}(\bx)e_{\lambda}(\by)=
            \sum_{d=0}^{\infty} \sum_{\substack{\lambda \vdash d\\\lambda_1 \leq N} } e_\lambda(\by)W^{\rm f}_\lambda(\bx)
       \end{split}
    \end{equation}
    Setting $\sum\limits_{i=1}^N y_i=0$ leads to shift invariance in $x$ and restricts the sum to $\lambda \in \mathcal{N}_N(d)$.
\end{proof}

\subsection{Clustering}
The clustering property was studied in relation to wavefunctions of the quantum Hall effect. It concerns vanishing properties under a $k+1$-fold coincidence of coordinates or $k+1$-\emph{clustering}.
\\

Consider, following \cite{feigin2002differential}, the space $ F_{N}^{(k)}$ of symmetric functions that vanish whenever $k+1$ variables coincide:
\begin{equation}
    f(z,\ldots z,x_{k+2},\ldots, x_N)=0 \quad \forall z
\end{equation}
The work \cite{feigin2002differential}  establishes a relation between Jack polynomials for certain admissible partitions and negative rational $\beta$ and this coincidence ideal. Let $J^{(\beta)}_{\lambda}$ denote a Jack polynomial of $N$ variables in monic normalization \footnote{Here $\beta=\alpha^{-1}$ compared to the $\alpha$ convention.}. Let $k \geq 1, r \geq 2$ and $r-1$ and $k+1$ be coprime. Set 
\begin{equation}
    \beta(r,k)= -\dfrac{r-1}{k+1}\,.
\end{equation}
A partition is called $(k,r,N)$ admissible if:
\begin{equation}
    \lambda_i-\lambda_{i+k} \geq r \,,\quad  1\leq i \leq N-k
\end{equation}
Typically, Jack polynomials will be singular at negative rational $\beta$. However, for $(k,r,N)$ admissible partitions, Jack polynomials of $N$ variables are finite at $\beta=\beta(r,k)$, given the coprimality condition above \cite{feigin2002differential}. It turns out that such Jack polynomials vanish under $k+1$-clustering. One of the results of \cite{feigin2002differential} is the following theorem.
\begin{theorem}[Feigin–Jimbo–Miwa–Mukhin 
\cite{feigin2002differential}] Taking $r=2$ in the construction above produces a basis for $k+1$-clustering polynomials
\begin{equation}
    F_{N}^{(k)} \simeq \operatorname{span}\left\{ J^{\beta(2,k)}_{\lambda}(\bx)  \, : \,  \lambda \text{ is } (k,2,N)\text{-admissible}\right\}
\end{equation}
\end{theorem}
The decay basis constructed in sec.\ref{sec::decay basis} behaves well under clustering, as determined by the following proposition.
\begin{proposition}
    For a nuclear partition $\lambda $ of $0<d\leq N$ and length $l=\ell(\lambda)$, the decay basis behaves as:
    \begin{equation}
  \begin{split}
    W_{\lambda}^{\rm d}(x_1,\ldots ,x_l\underbrace{,z,\ldots z}_{N-l}) = \dfrac{s_{\lambda}(x_1-z,\ldots, x_l-z)}{\sfC_{\lambda}(N)}=
    \\
    =\dfrac{\prod_{i=1}^l(x_i-z)^{\lambda_l}}{\sfC_{\lambda}(N)} s_{[\lambda_1-\lambda_l,\ldots \lambda_{l}-\lambda_l]}(x_1-z,\ldots ,x_l-z)
  \end{split}
  \end{equation}
\end{proposition}
Generally the admissible Jack polynomials are not shift invariant. The condition $\nabla_N J_{\lambda}=0$ is called the highest weight condition and is important in the context of the QHE. The conditions on $\lambda$ that give highest weight polynomials were found in \cite{bernevig2008generalized}. Therefore our decay polynomials form a basis in the space of highest weight polynomials with a fixed vanishing property.
\begin{lemma}
For $d\leq N$, let $F_{N,d}^{(k)} \subset F_{N}^{(k)}$ denote the degree $d$ subspace. Then:
    \begin{equation}
         \operatorname{span}\left\{W_{\lambda}^{\rm d}(\bx)\, :  \, \lambda \in \mathcal{N}(d), \ell(\lambda)\geq N-k  \right\}=\bar{\Lambda}_N \,\cap F_{N,d}^{(k)}
    \end{equation}
\end{lemma}
\begin{proof}
    For a shift invariant symmetric polynomial $g(\bx) \in \bar{\Lambda}_N$, the vanishing condition is equivalent to:
 
        \begin{equation}
    g(x_{1},\ldots x_{N-k-1},\underbrace{0,\ldots ,0}_{k+1})=0 \quad
    \end{equation}
    Under such specialization, Schur polynomials $s_{\mu}(\bx)$ vanish if $\ell(\mu) \geq N-k$, while the other Schur polynomials remain linearly independent as polynomials in the remaining $N-k-1$ variables. Therefore 
        $g(\bx) \in \bar{\Lambda}_N \cap F_{N,d}^{(k)} $ if and only if its Schur coefficients for $s_{\mu}$ with $\ell(\mu)<N-k$ vanish.
        \\
        
    First, we prove that decay polynomials with $\ell(\lambda) \ge N-k$ satisfy this condition. The Schur functions in the expansion of the decay element $\Wd_\lambda(\bx) $ all satisfy $\ell(\mu) \geq \ell(\lambda)$. Therefore if $\ell(\lambda) \geq N-k$, then the condition above is satisfied. 
        \\
        Conversely, take  $ g(\bx)\in \bar{\Lambda}_N \cap F_{N,d}^{(k)} $. Because $g(\bx)$ is shift invariant, expand it into the decay basis:
        \begin{equation}
            g(\bx)=\sum_{\mu \in \mathcal{N}(d)} g_\mu \Wd_\mu(\bx)
        \end{equation}
        Since the only nuclear partition in the Schur expansion of $\Wd_\lambda(\bx)$ is $\lambda$ itself, we get $g_\mu=\sfC_{\mu}(N)(s_\mu,g)$. Adding that $ g(\bx)\in F_{N,d}^k$ we get that :
        \begin{equation}
            g(\bx)=\sum_{\substack{\mu \in \mathcal{N}(d)
            \\
            \ell(\mu) \geq N-k}} g_\mu \Wd_\mu(\bx)
        \end{equation}
        This proves that decay polynomials with these restrictions indeed form a basis.
\end{proof}
The highest weight Jack polynomials and the decay basis lie in the same vector space. Therefore one could expect some relation between them. It turns out that in general the expansion of these rational Jacks in terms of the decay elements is complicated. Still, in some examples a simple relation can be found. 
\begin{corollary}
       For partitions  $\lambda=[2^{l}]$ we have the equality
\begin{equation}\label{eq:decaytoJack}
    \sfC_{[2^l]}(N)W_{[2^l ]}^{\rm d}(\bx)= J^{-(N-l+1)^{-1}}_{[2^l]}(\bfx) \,,\quad N\geq 2l
\end{equation}
where $r=2, k=N-l \geq l$.
\end{corollary}
\begin{proof}
    With this choice of parameters the partition $[2^l]$ is admissible. Moreover it is the unique admissible partition of degree $2l$ and there are no admissible partitions of lower degree. This means that $\dim F_{N,2l}^{(N-l)}=1$. Therefore, by the Lemma above, the respective decay and the Jack polynomials should be proportional. The coefficient $  \sfC_{[2^l]}(N)$ makes the leading monomial coefficient on the left of \eqref{eq:decaytoJack} to be one, hence establishing the equality. The Jack polynomial must also be highest weight, since the $\nabla_N$ operator preserves the clustering condition. Therefore $\nabla_N J^{-(N-l+1)^{-1}}_{[2^l]}(\bfx)$ must be in $F_{N,2l-1}^{(N-l)}$, which is empty. 
\end{proof}

\subsection{Antipodal self duality}
In ref.  \cite{Dixon:2026fuz} it was shown that shift invariant symmetric polynomials correspond to polynomials in so-called ladder integrals which are invariant under the antipode map of the Hopf algebra of multiple polylogarithms. 
\\

Let:
\begin{equation}
    \mathsf{f}_p=\sum_{j=p}^{2p} \frac{(p-1)! j! \, {\mathcal L}^{2p-j}}{(j-p)! (2p-j)!} \, \left[ {\rm Li}_j(z) - {\rm Li}_j(\bar{z}) \right] \,,
\end{equation}
where ${\mathcal L} \equiv -\log(z\bar{z})$, and ${\rm Li}_j(z) = \sum\limits_{k=1}^\infty \frac{z^k}{k^j}$ is the classical polylogarithm of weight $j$. These are the so-called ladder integrals (we use a different font to distinguish them from the forgotten symmetric functions). The coordinates $z,\bar{z}$ are the complex cross ratios of four points that parameterize conformal four point functions.
\\

There exists a map from the space  $\Lambda_m$ of symmetric polynomials in $m$ variables to the polynomial ring generated by $\mathsf{f}_p$ denoted as \cite{Basso:2021omx}:
\begin{equation}
    \ev{\ldots } :  \Lambda_m \rightarrow\mathbb{C}[\mathsf{f}_1,\mathsf{f_2},\ldots ]
\end{equation} For $F(\bx) \in \Lambda_m$
\begin{equation}\label{eq:fishnetEV}
    \ev{F(\bx)}=\frac{(-1)^{m(m-1)/2}}{m!}\int_{0}^{\infty}  \prod_{i=1}^m \mathrm{d}x_i\, w_{\rm F}(x_i)\, \prod_{i<j}(x_i-x_j)^2 \, F(\bfx)\,, 
\end{equation}
where the measure is given by:
\begin{equation}
   w_{\rm F}(x) = \left( \dfrac{z-\zb}{2\sqrt{z\zb}} \right)\cdot \dfrac{1}{\cosh\sqrt{x+\sigma^2}+\cosh\varphi} \,.
\end{equation}
and $\sigma, \varphi$ are defined as: 
\begin{equation}  
z= -\exp(\sigma+\varphi)\,, \qquad \zb= -\exp(\sigma-\varphi)\,.
\end{equation}
Define an automorphism $\hat{S}$ of the ring $\mathbb{C}[L,\mathsf{f}_1,\mathsf{f_2},\ldots ]$  as:
\begin{equation}
   \hat{S}(\mathsf{f}_p )=\sum_{k=0}^{p-1} \binom{p-1}{k} L^{k}\mathsf{f}_{p-k} 
\end{equation}
where $L=\log z \log \bar{z}$. The map extends to $L$ by requiring that $\hat{S}(L)=-L$. The main result of \cite{Dixon:2026fuz} is that it can be expressed in terms of the $L$-dependent shift of the symmetric function:
\begin{equation}
    \hat{S}\ev{F(\bx)}=\ev{F(\bx+L{\bf 1})}
\end{equation}
The ladder integrals are linear combinations of single valued polylogarithms. The map $\hat{S}$ descends from the antipode map on the Hopf algebra of multiple polylogarithms \cite{Goncharov2001,Brown:2015ylf,Duhr:2012fh}. It is then composed with a kinematic map, which is explained in \cite{Dixon:2025zwj,Dixon:2026fuz}.
The correspondence induced by the integral \eqref{eq:fishnetEV} implies invariants of $\hat{S}$, called antipodal invariants, can be obtained from shift invariant symmetric polynomials. In \cite{Dixon:2026fuz} the decay basis was used to find antipodal invariants $Q_{\lambda}^{(m)}=\ev{W^{\rm decay}_{\lambda}}$, which are defined  for $|\lambda| \leq m$. The map  \eqref{eq:fishnetEV} has an important property which motivates the study of Schur expansions. Namely, by virtue of Andréief’s identity \cite{andreief1883note,forrester2019meet} one can explicitly calculate its value on Schur polynomials:
\begin{equation}
    \ev{s_{\lambda}(\bx)}=\det {\sf f}_{\lambda_i+j-i+m}
\end{equation}
The determinants are denoted as:
\begin{equation}
N_{\lambda}^{(m)}=\det {\sf f}_{\lambda_i+j-i+m}
\end{equation}
Therefore in order to establish explicit expressions for antipodal invariants one needs to compute the Schur expansion of shift invariant symmetric polynomials. As a corollary of the results of this paper we find three more explicit formulas for antipodal invariants. Crucially, the one provided by the fusion basis works for all $m$ and $\lambda \in \mathcal{N}_m$.
\begin{lemma}
    The following ladder polynomials are invariant under the map $\hat{S}$:
    \begin{equation}
        \begin{split}
Q^{(m),i}_{\lambda} = \ev{W_{\lambda}^{i}(\bx) } \,,\quad i={\rm f,F,d,p}
        \end{split}
    \end{equation}
    where $Q^{(m),\rm f}_{\lambda},Q^{(m),\rm F}_{\lambda}$ are defined for any $\lambda \in \mathcal{N}_m$, and $Q^{(m),\rm d}_{\lambda}$, $Q^{(m),\rm p}_{\lambda}$  are defined only for $\lambda \in \mathcal{N}_m$ with $|\lambda|\leq m$.
\end{lemma}
The invariants $Q^{(m),\rm f}_{\lambda}$ or $Q^{(m),\rm F}_{\lambda}$ solve the problem of resolving the $m$ singularities, raised in \cite{Dixon:2026fuz} . In particular, the degree 4 example in that paper asserts that the combination $Q^{(m),\rm d}_{[4]}+Q^{(m),\rm d}_{[2,2]}$ is well defined at $m=2,3$. In notations of this section it is nothing but:
\begin{equation}
   Q_{[2,2]}^{(m),\rm f}= Q^{(m),\rm d}_{[4]}+Q^{(m),\rm d}_{[2,2]}
\end{equation}

\section{Acknowledgments}

V.M. thanks L. Dixon, C. Duhr, and F. Fernandes for collaboration on earlier work, which inspired further investigations into shift invariant symmetric polynomials. We thank A. Morozov, A. Mironov, A. Sleptsov, A. Popolitov and P. Alexandersson for interesting discussions. The authors used GPT-6 Astra for help with some of the proofs and for polishing the exposition and spellchecking. When we initially constructed the preimage of the fusion basis, we (ironically) didn't recognize the forgotten symmetric functions straight away. GPT-6 was very helpful at that point. All the LLM responses were manually verified, and the final text was written by the authors. This work is supported by the RSF grant 26-12-00191. 

\bibliographystyle{utphys}
\bibliography{references}

\appendix

\section{Properties of forgotten symmetric functions}
\label{sec::append A}
It is instructive to study properties of forgotten symmetric functions via the usual monomial symmetric functions
\begin{equation}
    m_{\lambda}(\bx) = \sum_{\text{symm}} x_{1}^{\lambda_1} x_2^{\lambda_2} \ldots x_{N}^{\lambda_N}
\end{equation}
For our purposes it is useful to proceed to augmented monomials \cite{Macdonald:book} $\tilde{m}_{\lambda} =  a_{\lambda} \cdot m_{\lambda}$. Recall that  $a_{\lambda}=\prod\limits_{j \geq 1} m_j(\lam)!$. Note that $m_{i}(\lambda)$ denotes the multiplicity of $i$ in $\lambda$. The multiplication rule for augmented monomials \cite{Merca2015} is:
\begin{equation}
    \tilde{m}_{[k]} \tilde{m}_{\lambda} = \sum_{i = 1}^{\ell(\lambda) + 1} \tilde{m}_{: \lambda + k \delta_i :}
\end{equation}
where $ : X : $ denotes the partition obtained by ordering $X$ in decreasing order, such that the result is a Young diagram. $\lambda + k \delta_i$ denotes the sequence obtained by adding $k$ to the $i$-th component of $\lambda$. For example,
\begin{equation}
    \tilde{m}_{[1]} \tilde{m}_{[1,1]} = \tilde{m}_{: [ 2,1] :} + \tilde{m}_{: [1,2] :} + \tilde{m}_{:[1,1,1]:} = 2 \tilde{m}_{[2,1]} + \tilde{m}_{[1,1,1]}
\end{equation}
This multiplication rule can be used to efficiently compute these polynomials by recursion as functions of $p_k$-variables. In other words,
\begin{equation}
    \tilde{m}_{[\lambda_1, \lambda_2, \ldots, \lambda_{\ell(\lambda)}, k]} = \tilde{m}_{[k]} \tilde{m}_{\lambda} - \sum_{i = 1}^{\ell(\lambda)} \tilde{m}_{: \lambda + k \delta_i :}
\end{equation}
On the r.h.s. diagrams of length not greater than $\ell(\lambda)$ appear, however on the l.h.s. one has a diagram of length $\ell(\lambda) + 1$. The initial data for the recursion is given by $\tilde{m}_{[k]} = p_k$.
As functions of $p_k$-variables, the rescaled monomial functions are related to the forgotten polynomials by the involution $\omega$ \cite{Macdonald:book}. Therefore,
    \begin{equation}
    f_{\lambda} (p) =  \frac{\tilde{m}_{\lambda} \left((-1)^{k+1}p_k\right)}{a_\lambda}  = m_{\lambda} \left((-1)^{k+1}p_k\right)
\end{equation}
\begin{proposition}
    Forgotten symmetric functions as functions of power sums possess the following generating function
    \begin{equation}
    \label{gen funct for forgotten}
        \exp\left( \sum_{\substack{\lambda \in \mathcal{P}
        \\
        \lambda \neq \varnothing}}(-1)^{|\lam| + \ell(\lam)} v_{\lambda} \, p_{|\lambda|} \cdot t_{\lambda}  \right) = \sum_{\mu} f_{\mu}(p) \cdot t_{\mu}
    \end{equation}
    where $t_{\lambda} = \prod\limits_i t_{\lambda_i}$ is a monomial in auxiliary variables $t_k$ and  $v_{\lambda} = \frac{(\ell(\lam) - 1)!}{\prod\limits_{i \geq 1} m_i(\lam)!}$.
\end{proposition}
\begin{proof}
First we note 
\begin{equation}
    \prod_{i=1}^{N} \left( 1 + \sum_{n = 1}^{\infty} (x_i)^n \cdot t_n \right) = \sum_{\mu} m_{\mu} (\bfx) \cdot t_{\mu} 
\end{equation}
On the other hand we can treat the l.h.s. via the exponent
\begin{equation}
    \prod_{i=1}^{N} \left( 1 + \sum_{n = 1}^{\infty} (x_i)^n \cdot t_n \right) = \exp \left( \sum_{i=1}^{N} F(x_i) \right)
\end{equation}
where
\begin{equation}
    F(x) = \log \left( 1 + \sum_{n = 1}^{\infty} x^n \cdot t_n \right) = \sum_{k = 1}^{\infty} x^k \sum_{s = 1}^{\infty} \frac{(-1)^{s-1}}{s} \sum_{\substack{n_1 + n_2 + \ldots + n_s = k \\ n_i \geq 1}} t_{n_1} t_{n_2} \ldots t_{n_s}
\end{equation}
Then one can reformulate the sum over compositions as a sum over partitions and obtain the generating function
\begin{equation}
     \exp\left( \sum_{\substack{\lambda \in \mathcal{P}
        \\
        \lambda \neq \varnothing}} (-1)^{\ell(\lam) - 1} v_{\lambda} \, p_{|\lambda|}(\bfx) \cdot t_{\lambda}  \right) = \sum_{\mu} m_{\mu} (\bfx) \cdot t_{\mu} 
\end{equation}
Forgotten symmetric functions then possess a generating function due to the involution. In particular, under the involution power sums transform as $\omega(p_k) = (-1)^{k+1} p_k$.
\end{proof}
From the above generating function directly follows the following interesting interpretation:
\begin{proposition}
    Forgotten symmetric functions obey the following fusion formula
    \begin{equation}
        f_{\mu} (p) = (-1)^{|\mu| + \ell(\mu)} \sum_{s = 1}^{\infty} \sum_{\substack{: \lambda^{(1)} \cup \lambda^{(2)} \ldots \cup \lambda^{(s)} : = \mu \\ \ell(\lam^{(i)}) \geq 1 }} \frac{\prod_{i = 1}^{s} v_{\lam^{(i)}} p_{|\lam^{(i)}|} }{s!} 
    \end{equation}
\end{proposition}
Note that in the above formula $\lam^{(i)}$ are partitions, in contrast to $\lam_i$ which denotes the $i$-th part of the partition. The term \textit{fusion} comes from the following picture: in each term in the sum, for each subset $\lambda^{(i)}$ all its components ``fuse'' into one, namely, $p_{\lam^{(i)}} = \prod\limits_{j \in \lambda^{(i)}}p_{j}$ goes to $p_{|\lambda^{(i)}|}$.

Formula \eqref{gen funct for forgotten} shows that $\epsilon_\lambda f_\lambda$  has non-negative coefficients in the power sum basis. Also due to the property $v_{r \cdot \lambda} = v_{\lambda}$ of the coefficients, the forgotten functions satisfy the following plethystic identity:
\begin{equation}
    f_{r \cdot \lambda  } (p_a) = (-1)^{|\lam|(r-1)} f_{\lam} (p_{r\cdot a})
\end{equation}
On the other hand, this property is a simple consequence of the plethystic property of monomial symmetric functions, which can be seen from the definition. We provide a small set of explicit examples
\begin{align}
    \begin{aligned}
        f_{[1^{k}]} &= s_{[k]} \\
        f_{[2,1]} &= -(p_1 p_2+p_3) \\
        f_{[3,1]} &= p_1 p_3+p_4 \\
        f_{[4,1]} &= -(p_1 p_4+p_5) \\
        f_{[3,2]} &= -(p_2 p_3+p_5) \\
        f_{[2,1,1}] &= -\frac{1}{2} \left(p_1^2 p_2+p_2^2+2 p_1 p_3+2 p_4\right) \\
        f_{[3,1,1]} &= \frac{1}{2} \left(p_1^2 p_3+p_2 p_3+2 p_1 p_4+2 p_5\right) \\
        f_{[2,2,1}] &=\frac{1}{2} \left(p_1 p_2^2+2 p_2 p_3+p_1 p_4+2 p_5\right) \\
        f_{[2,1,1,1]} &= -\frac{1}{6} \left(p_1^3 p_2+3 p_1 p_2^2+3 p_1^2 p_3+5 p_2 p_3+6 p_1 p_4+6 p_5\right)
    \end{aligned}
\end{align}

\end{document}